\documentclass[aop,MSNbibl,seceqn,dvips]{arximspdf}

\doi{10.1214/12-AOP793} 
\volume{41}
\issue{5}
\pubyear{2013}
\firstpage{3345}
\lastpage{3391}

\makeatletter
\newcommand{\rrvert}{\vert}
\newcommand{\llvert}{\vert}

\newtheorem{lemma}{Lemma}
\newtheorem{theorem}{Theorem}
\newtheorem{corollary}{Corollary}
\newtheorem{proposition}{Proposition}

\newproclaim{definition}{Definition}

\newcommand{\s}{\bolds{\sigma}}
\newcommand{\T}{\operatorname{th}}
\newcommand{\C}{\operatorname{ch}}
\newcommand{\sh}{\operatorname{sh}}
\newcommand{\w}{\omega}

\makeatother

\begin{document}
\begin{frontmatter}

\title{Disorder chaos in the Sherrington--Kirkpatrick model with
external field}
\runtitle{Disorder chaos in the SK model}

\begin{aug}
\author[A]{\fnms{Wei-Kuo} \snm{Chen}\corref{}\ead[label=e1]{weikuoc@uci.edu}}
\runauthor{W.-K. Chen}
\affiliation{University of California, Irvine}
\address[A]{Department of Mathematics\\
340 Rowland Hall\\
University of California, Irvine\\
Irvine, California 92697-3875\\
USA\\
\printead{e1}} 
\end{aug}

\received{\smonth{9} \syear{2011}}
\revised{\smonth{7} \syear{2012}}

%
\begin{abstract}
We consider a spin system obtained by coupling two distinct
Sherring\-ton--Kirkpatrick (SK) models with the same temperature and
external field whose Hamiltonians are correlated. The disorder chaos
conjecture for the SK model states that the overlap under the
corresponding Gibbs measure is essentially concentrated at a single
value. In the absence of external field, this statement was first
confirmed by Chatterjee [Disorder chaos and multiple valleys in spin
glasses (2009) Preprint]. In the present paper, using Guerra's replica
symmetry breaking bound, we prove that the SK model is also chaotic in
the presence of the external field and the position of the overlap is
determined by an equation related to Guerra's bound and the Parisi
measure.
\end{abstract}

%
\begin{keyword}[class=AMS]
\kwd{60K35}
\kwd{82B44}
\end{keyword}
\begin{keyword}
\kwd{Disorder chaos}
\kwd{Guerra's replica symmetry breaking bound}
\kwd{Parisi formula}
\kwd{Parisi measure}
\kwd{Sherrington--Kirkpatrick model}
\end{keyword}

\end{frontmatter}

\section{Introduction and main results}
The phenomenon of chaos arose from the discovery that in some models, a
slight perturbation on the parameters such as the temperature, external
field or disorder will result in a dramatic change to the system. In
this paper, we will be concerned with the Sherrington--Kirkpatrick (SK)
model~\cite{SK75} and study its chaotic property mainly due to the
change of the disorder.
Let us begin by recalling the definition of the SK model and the
formulation of the Parisi formula.
Suppose that $\xi\dvtx \Bbb{R}\rightarrow\Bbb{R}$ is a convex function satisfying
$\xi(x)=\xi(-x)$, $\xi''(x)>0$ if $x\neq0$, and $\xi^{(3)}\geq0$
if $x>0$.
For each $N$, we consider a centered Gaussian process $H=H_N$ indexed
by the configuration space $\Sigma_N= \{-1,+1 \}^N$ with
covariance
\[
EH_N\bigl(\s^1\bigr)H_N\bigl(
\s^2\bigr)=N\xi(R_{1,2})
\]
for
$\s^1=(\sigma_1^1,\ldots,\sigma_N^1),\s^2=(\sigma_1^2,\ldots,\sigma_N^2)\in\Sigma_N$,
where
\[
R_{1,2}=R_{1,2}\bigl(\s^1,\s^2\bigr)=
\frac{1}{N}\sum_{i\leq N}\sigma_i^1
\sigma_i^2
\]
is called the overlap of the configurations $\s^1$ and $\s^2$. Let
$h$ be a random variable and $(h_i)_{i\leq N}$ be i.i.d. copies of $h$.
Then the SK model with external field $h$ possesses the Hamiltonian
\[
-H(\s)+\sum_{i\leq N}h_i\sigma_i
\]
for $\s=(\sigma_1,\sigma_2,\ldots,\sigma_N)\in\Sigma_N$ and its
Gibbs measure is defined as
\[
G_N(\s)=\frac{1}{Z_N}\exp\biggl(-H(\s)+\sum
_{i\leq N}h_i\sigma_i \biggr),
\]
where $Z_N$ is a normalizing factor, called the partition function. Let
us also define
\[
p_N=\frac{1}{N}E\log Z_N=\frac{1}{N}E\log
\sum_{\s\in\Sigma
_N}\exp\biggl(-H(\s)+\sum
_{i\leq N}h_i\sigma_i \biggr).
\]
This quantity is usually called the free energy for the SK model in
physics and its thermodynamic limit $\lim_{N\rightarrow\infty}p_N$
can be computed by the Parisi formula described below.

Consider an integer $k\geq0$ and numbers
%
\begin{eqnarray}
\label{INTeq2}
&&\mathbf{m}\mbox{: } m_0=0\leq m_1\leq\cdots\leq
m_k\leq m_{k+1}=1,
\nonumber\\[-8pt]\\[-8pt]
&&\mathbf{q}\mbox{: } q_0=0\leq q_1\leq\cdots\leq q_{k+1}
\leq q_{k+2}=1.
\nonumber
\end{eqnarray}
It helps to think of the triplet $k,\mathbf{m},\mathbf{q}$ as a
probability measure $\mu$ on $ [0,1 ]$ that has all its
mass concentrated at a finite
number of points $q_1,\ldots,q_{k+1}$ and $\mu( [0,q_p
])=m_p$ for $1\leq p\leq k+1$.
Let $z_0,\ldots,z_{k+1}$ be independent Gaussian r.v.'s with
$Ez_p^2=\xi'(q_{p+1})-\xi'(q_p)$ for $0\leq p\leq k+1$.
Starting with
\[
X_{k+2}=\log\cosh\biggl(h+\sum_{0\leq p\leq k+1}z_p
\biggr),
\]
we define by decreasing induction for $1\leq p\leq k+1$,
\[
X_p=\frac{1}{m_p}\log E_p\exp m_pX_{p+1},
\]
where $E_p$ means the expectation on the r.v.'s
$z_p,z_{p+1},\ldots,z_{k+1}$. If $m_p=0$ for some $p$, we define
$X_p=E_{p}X_{p+1}$. Finally, we define $X_0=EX_1$. Set
%
\begin{equation}
\label{INTeq5} \mathcal{P}_k(\mathbf{m},\mathbf{q})=
\log2+X_0-\frac{1}{2}\sum_{p=1}^{k+1}m_p
\bigl(\theta(q_{p+1})-\theta(q_{p})\bigr),
\end{equation}
where $\theta(x)=x\xi'(x)-\xi(x)$. This quantity is the famous
Guerra replica symmetry breaking bound of the $k$th level~\cite{G03}
that yields a fundamental inequality, for every $k,\mathbf{m},\mathbf{q}$,
%
\begin{equation}
\label{addeq7} p_N\leq\mathcal{P}_k(\mathbf{m},
\mathbf{q}).
\end{equation}
Let us define the Parisi functional on the space of all probability
measures on $[0,1]$ consisting of only a finite number of point masses
by $\mathcal{P}(\xi,h,\mu)=\mathcal{P}_k(\mathbf{m},\mathbf{q})$
if $\mu$ corresponds to $(k,\mathbf{m},\mathbf{q})$. We define
$\mathcal{P}(\xi,h)=\inf_{k,\mathbf{m},\mathbf{q}}\mathcal
{P}_k(\mathbf{m},\mathbf{q})$,
where the infimum is over all choices of $(k,\mathbf{m},\mathbf{q})$
as above.
Then the Parisi formula says that
\[
\lim_{N\rightarrow\infty}p_N=\mathcal{P}(\xi,h).
\]
This formula was first rigorously proven in Talagrand~\cite{Talag062}.
It is well known~\cite{G03} that the Parisi functional is Lipschitz
continuous with respect to the metric $d(\mu,\mu')=\int_0^1|\mu
([0,q])-\mu'([0,q])|\,dq$. Thus, it can be extended continuously to the
space of all probability measures defined on $[0,1]$ and is denoted
again by $\mathcal{P}(\xi,h,\cdot)$. Then clearly $\lim_{N\rightarrow
\infty}p_N=\mathcal{P}(\xi,h)=\min\mathcal{P}(\xi,h,\mu)$, where the
minimum is taken over all probability measures defined on $[0,1]$. Any
measure that achieves the minimum is called a Parisi measure.
Heuristically, one may think of the Parisi measure as the limiting
distribution of the overlap.

We are now ready to formulate the disorder chaos problem in the SK
model. Let $0\leq t\leq1$.
Suppose that $H^1=H_N^1$ and $H^2=H_N^2$ are two centered Gaussian
processes having the same distribution as $H$
and they are correlated in the following way,
%
\begin{equation}
\label{INTeq1} EH^1\bigl(\s^1\bigr)H^2\bigl(
\s^2\bigr)=Nt\xi(R_{1,2}).
\end{equation}
That is, we allow a portion $1-t$ of independence between two systems.
Consider the coupled Hamiltonian
\[
-H^1\bigl(\s^1\bigr)-H^2\bigl(
\s^2\bigr)+\sum_{i\leq N}h_i\bigl(
\sigma_i^1+\sigma_i^2\bigr)
\]
on $\Sigma_N^2$. Proceeding as before, we define its Gibbs measure by
\[
G_N' \bigl(\s^1,\s^2 \bigr)=
\frac{1}{Z_{N}'}\exp\biggl(-H^1\bigl(\s^1
\bigr)-H^2\bigl(\s^2\bigr)+\sum
_{i\leq N}h_i\bigl(\sigma_i^1+
\sigma_i^2\bigr) \biggr),
\]
where the normalizing factor $Z_N'$ is the partition function of this
model. As we have already mentioned in the beginning of this section,
the chaos phenomenon is concerned with the instability occurring in
some spin glass models due to the change of some external parameters.
In the SK model, one very basic way to measure such instability mainly
due to the change of the disorder, or, briefly, chaos in disorder, is
to study the behavior of the overlap. A typical statement one is
looking for in this case is that if $0<t<1$, the overlap takes
essentially only one value under $G_N'$. This is quite different from
the typical lack of self-averaging property of the overlap in the low
temperature phase when $t=1$. The phenomenon of chaos itself was first
conjectured by Fisher and Huse~\cite{FH86}. Early discussion on the
disorder chaos for the SK model can be found in~\cite{BM87} and \cite
{MBK82}. For further references in the physics literature, one may
refer to~\cite{KK07}.
However, the mathematically rigorous results have appeared only lately.
In the absence of the external field, Chatterjee~\cite{Cha09} recently
proved chaos in disorder and discovered that the overlap is
concentrated at $0$.

In the present work, we aim to prove that the disorder chaos also holds
in the presence of the external field, that is, $Eh^2\neq0$.
Moreover, we find that when there is chaos, the position of the overlap
can be described by an equation, which is
related to the Parisi measure and can be formulated as follows. Suppose
that $\mu$ is a Parisi measure. Recall that $\mu$ minimizes the
Parisi functional. We can approximate $\mu$ weakly by a sequence of
$\varepsilon_n$-stationary measures $(\mu_n)$ satisfying $\mathcal
{P}(\xi,h,\mu_n)\rightarrow\mathcal{P}(\xi,h)$. Here, by
$\varepsilon_n$-stationarity, it means that the measure $\mu_n$
minimizes the $k$th level Guerra replica symmetry breaking bound for
some $k$ depending on $n$ and $\mathcal{P}(\xi,h,\mu_n)<\mathcal
{P}(\xi,h)+\varepsilon_n$, where $\varepsilon_n\downarrow0$ (see
Definition~\ref{sec3def2} below). This approximation is for technical
purposes that have played a crucial role in Talagrand's proof on the
Parisi formula~\cite{Talag062} and will also be of great importance in
our argument. For a given $(k,\mathbf{m},\mathbf{q})$ corresponding
to $\mu$, recall the definition of $X_0$ from (\ref{INTeq5}).
A~very nice and useful fact about this quantity is that it can be
computed as $E\Phi(h,0)$, where $\Phi\dvtx \mathbb{R}\times
[0,1]\rightarrow\mathbb{R}$ is the solution to the following PDE,
%
\begin{equation}
\label{INTeq3}\quad \frac{\partial\Phi}{\partial q}=-\frac{\xi''(q)}{2}
\biggl(\frac
{\partial^2\Phi}{\partial x^2}+\mu
\bigl( [0,q ]\bigr) \biggl(\frac
{\partial\Phi}{\partial
x} \biggr)^2 \biggr)\qquad
\forall(x,q)\in\mathbb{R}\times[0,1]
\end{equation}
with $\Phi(x,1)=\log\cosh x$. For each $n$, let $\Phi_n$ be the PDE
solution (\ref{INTeq3}) corresponding to $\mu_n$.
From~\cite{Talag06}, we know that $(\Phi_n)$ converges uniformly and
we denote its limit by $\Phi$. Moreover,~\cite{Talag06} yields that
the first partial derivative of $\Phi$ with respect to $x$ exists.
From this, for each fixed $0< v<1$, we define
%
\begin{equation}
\label{INTeq4} \varphi_v(u,t)=E\frac{\partial\Phi}{\partial x}(h+
\chi_1,v)\,\frac
{\partial\Phi}{\partial x}(h+\chi_2,v)-u
\end{equation}
for all $0\leq u\leq v$ and $0\leq t\leq1$, where $\chi_1$ and $\chi_2$
are jointly Gaussian with $E\chi_1^2=E\chi_2^2=\xi'(v)$ and
$E\chi_1\chi_2=t\xi'(u)$ independent of $h$. The motivation of
$\varphi_v$ comes from the Guerra replica symmetry breaking bound for
the coupled free energy that will be explained in great detail in
Section~\ref{sec5} below. An important fact about the Parisi measure
$\mu$ in the case of $Eh^2\neq0$ is that the smallest value $c$ of
its support is positive. This is called the positivity of the overlap
(see Chapter 14~\cite{Talag11}). When $v=c$ and $0\leq t<1$, we are
able to determine the number of the solutions of $\varphi_c(\cdot,t)=0$.

\begin{proposition}\label{ME0cprop1}
For each $0\leq t<1$, there exists a unique $u_t$ in $ [0,c
]$ such that $\varphi_c(u_t,t)=0$.
Moreover, $\varphi_c(c,t)<0$ for $0\leq t<1$ and $\varphi_c(c,1)=0$.
\end{proposition}

Now, the quantitative result of the disorder chaos in the SK model is
stated as follows.

\begin{theorem}\label{INTthm1}
Suppose that $0<t<1$ and $Eh^2>0$. Then the SK model has disorder
chaos, namely, for any $\varepsilon>0$, the following holds:
%
\begin{equation}
\label{INTthm1eq1} EG_N' \bigl( \bigl\{\bigl(
\s^1,\s^2\bigr)\dvtx |R_{1,2}-u_t|\geq
\varepsilon\bigr\} \bigr)\leq K\exp\biggl(-\frac{N}{K} \biggr),
\end{equation}
where $K$ is a constant depending on $t,\xi,h,\varepsilon$ and $\mu$.
\end{theorem}

A consequence of Theorem~\ref{INTthm1} is that even though we do
not know that the Parisi measure $\mu$ is unique,
the quantity $u_t$ is independent of the choice of $\mu$. However, the
convergence rate $K$ in (\ref{INTthm1eq1}) does depend on $\mu$.
In~\cite{Talag06} and~\cite{Talag07}, other types of chaos
problems in the SK model
are also proposed, such as chaos in temperature and chaos in external
field. Again, the rigorous results are still scarce.
Theorem~\ref{INTthm1} is the first result in chaos problems of any
kind in the SK model with the external field.
To the best of our knowledge, the only other two instances of chaos
problems in spin glasses are in the work of
Chatterjee~\cite{Cha09}, who proved chaos in disorder in the SK model
without the external field,
and in the work of Panchenko and Talagrand~\cite{PT07}, who
established chaos in the external field in the spherical SK model.

The approach of the present paper is motivated by Talagrand's proof on
the positivity of the overlap in the SK model; see Section 14.12
\cite{Talag11}. We also refer to a sketch of a possible proof for the
disorder chaos problem discussed in Research Problem 15.7.14 \cite
{Talag11}. However, it is by no means clear how to implement these
approaches properly that contain several technical issues and require
some new ideas. Here is our main result.

\begin{proposition}\label{INTprop1}
Let $0<t<1$ and $Eh^2>0$. For $\varepsilon>0$, there exists some
$\varepsilon^*>0$ such that
%
\begin{eqnarray}
\label{INTprop1eq1}\qquad p_{N,u}:\!&=&\frac{1}{N}E\log\sum
_{R_{1,2}=u}\exp\biggl(-H^1\bigl(\s^1
\bigr)-H^2\bigl(\s^2\bigr)+\sum
_{i\leq N}h_i\bigl(\sigma_i^1+
\sigma_i^2\bigr) \biggr)
\nonumber\\[-8pt]\\[-8pt]
&\leq&2\mathcal{P}(\xi,h)-\varepsilon^*
\nonumber
\end{eqnarray}
for all $u$ satisfying $|u-u_t|\geq\varepsilon$, where $\varepsilon^*$
is a constant depending on $t,\xi,h,\varepsilon$ and~$\mu$.
\end{proposition}

As an immediate consequence of the Gaussian concentration of measure
phenomenon (see Theorem 13.4.3 in~\cite{Talag11} and also the argument
for the positivity of the overlap on page 449 of~\cite{Talag11}),
Theorem~\ref{INTthm1}
follows from Proposition~\ref{INTprop1}. Let us continue by giving a
brief description\vadjust{\goodbreak} of how we proceed to prove Proposition~\ref{INTprop1}.
The approach for proving (\ref{INTprop1eq1}) is based on the
Guerra replica symmetry breaking bound
that was first used for the coupled system in~\cite{Talag062}. We
divide our discussion into three cases: $-1\leq u\leq0$, $0\leq u\leq
c'$, and $c'<u\leq1$,
where $c'$ satisfies $c'>c$ and is very close to $c$.
In the presence of the external field, we adapt a similar argument as
Talagrand's proof on the positivity of the overlap (see Section $14.12$
in~\cite{Talag11}) to conclude (\ref{INTprop1eq1}) for $-1\leq
u\leq0$.
In the case that $0\leq u\leq c'$, if there is chaos, the system should
exhibit ``high temperature behavior''
and $u_t$ should be determined by an equation related to the Parisi
measure, as is the case of the original SK model
in the high temperature regime; see Chapter 2 in~\cite{Talag101}. This
observation then leads to (\ref{INTprop1eq1}).
The most difficult part of our study is the case when $c'<u\leq1$. We
establish an iterative inequality,
which is very sensitive to the parameter $t$. From the construction of
the Parisi measure, we are able to find
parameters such that (\ref{INTprop1eq1}) holds even in the absence
of the external field.

The paper is organized as follows. Throughout the paper, we denote by
$E$ the expectation with respect to all randomness and we assume that
the external field $h$ satisfies $Eh^2>0$ and every Gaussian r.v. is
centered. In Section~\ref{sec1} we first give the formulation of an
extended version of Guerra's replica symmetry breaking bound
and explain why this is applicable to our study. We then continue to
carry out the core of the proof of Proposition~\ref{INTprop1}.
In Section~\ref{sec3} we state some results that help to control
Guerra's bound.
Most of their proofs can be found in~\cite{Talag11}.
Section~\ref{sec4} is devoted to proving (\ref{INTprop1eq1})
for $-1\leq u\leq0$ based on the same argument as
Section~14.12 in~\cite{Talag11}. In Section~\ref{sec5} we study how
Guerra's bound relates to the definition of $\varphi_v$ and
give the proof of Proposition~\ref{ME0cprop1}. Together they imply
(\ref{INTprop1eq1}) for $0\leq u\leq c'$.
Finally, we develop an iterative inequality and prove
(\ref{INTprop1eq1}) for $c'<u<1$ in Section~\ref{sec6}.


\section{Methodology}\label{sec1}

Let us first state an extension of the Guerra replica symmetry breaking bound.
Suppose that $-1\leq u\leq1$ and $\eta\in\{-1,+1\}$
satisfies $u=\eta|u|$.
For a given integer $\kappa\geq1$, we consider numbers
%
\begin{eqnarray}
\label{GBasseq1}
1&\leq&\tau\leq\kappa,\qquad\tau\in\Bbb{N},
\nonumber
\\
n_0&=&0\leq n_1\leq\cdots\leq n_{\kappa-1}\leq
n_\kappa=1,
\\
\rho_0&=&0\leq\rho_1\leq\cdots\leq\rho_\tau=|u|
\leq\rho_{\tau
+1}\leq\cdots\leq\rho_{\kappa+1}=1.
\nonumber
\end{eqnarray}
For $0\leq p\leq\kappa$, suppose that we are given independent pairs
of jointly Gaussian r.v.'s $(y_p^1,y_p^2)$ with
\[
E\bigl(y_p^1\bigr)^2=E\bigl(y_p^2
\bigr)^2=\xi'(\rho_{p+1})-\xi'(
\rho_p)
\]
such that
\[
Ey_p^1y_p^2=\eta t\bigl(
\xi'(\rho_{p+1})-\xi'(\rho_p)
\bigr) \qquad\mbox{if $0\leq p< \tau$}
\]
and
\[
\mbox{$y_p^1$ and $y_p^2$ are
independent if $\tau\leq p\leq\kappa$.}
\]
These r.v.'s are independent of $h$.
For our convenience, from now on, we set $\sh(x)=\sinh x$, $\C
(x)=\cosh x$, and $\T(x)=\tanh x$.
Let $\lambda$ be any real number. Starting with
\begin{eqnarray*}
Y_{\kappa+1}&=&\log\biggl(\C\biggl(h+\sum_{0\leq p\leq\kappa
}y_p^1
\biggr)\C\biggl(h+\sum_{0\leq p\leq\kappa}y_p^2
\biggr)\C\lambda
\\
&&{} +\sh\biggl(h+\sum_{0\leq p\leq\kappa
}y_p^1
\biggr)\sh\biggl(h+\sum_{0\leq p\leq\kappa}y_p^2
\biggr)\sh\lambda\biggr),
\end{eqnarray*}
we define by decreasing induction for $p\geq1$,
\[
Y_p=\frac{1}{n_p}\log E_p\exp
n_pY_{p+1},
\]
where $E_p$ denotes expectation in the r.v.'s $y_n^j$ for $n\geq p$. In
the case of $n_p=0$ for some $p$, we set
$Y_p=E_pY_{p+1}$. Finally, we define $Y_0=EY_1$.

\begin{theorem}\label{GBthm2}
We have
%
\begin{eqnarray}
\label{GBthm2eq1} p_{N,u}&\leq&2\log2+Y_0-\lambda u-(1+t)
\sum_{0\leq p<\tau}n_p\bigl(\theta(
\rho_{p+1})-\theta(\rho_p)\bigr)
\nonumber\\[-8pt]\\[-8pt]
&&{}-\sum_{\tau\leq p\leq\kappa}n_p\bigl(\theta(
\rho_{p+1})-\theta(\rho_p)\bigr).
\nonumber
\end{eqnarray}
\end{theorem}
Recalling Guerra's original bound (\ref{addeq7}),
(\ref{GBthm2eq1}) is a kind of two-dimensional extension.
Its proof is essentially the same as that of Proposition 14.12.4 \cite
{Talag11} and a more generalized version can be found in Section 15.7
\cite{Talag11}. One might have already observed that from the
definition of $p_{N,u}$ and (\ref{addeq7}), $p_{N,u}\leq2p_N\leq
2\mathcal{P}_k(\mathbf{m},\mathbf{q})$ for any $k,\mathbf
{m},\mathbf{q}$. Before we proceed to state our main results in this
section, let us illustrate that for any given $k,\mathbf{m},\mathbf
{q}$, we can find parameters (\ref{GBasseq1}) such that the
right-hand side of (\ref{GBthm2eq1}) is equal to $2\mathcal
{P}_k(\mathbf{m},\mathbf{q})$. This recovers the inequality
$p_{N,u}\leq2\mathcal{P}_k(\mathbf{m},\mathbf{q})$. To do this, let
$k,\mathbf{m},\mathbf{q}$ satisfy (\ref{INTeq2}) and $\tau$ with
$1\leq\tau
\leq k+2$ satisfying
\[
q_{\tau-1}\leq|u|\leq q_\tau.
\]
Without loss of generality, we may assume that $|u|$ is in the list of
$\mathbf{q}$. Indeed, we can always consider
a new triplet $k+1,\mathbf{m}',\mathbf{q}'$ obtained by inserting
$|u|$ into $\mathbf{q}$ and keeping $\mathbf{m}$ fixed in the
following way:
\begin{eqnarray*}
&&\mathbf{m'}\mbox{: $m_p'=m_p$
for $0\leq p\leq\tau-1$, $m_{p-1}$ if $p=\tau$,
and $m_{p-1}$ if $\tau+1\leq p\leq k+2$,}
\\
&&\mathbf{q'}\mbox{: $q_p'=q_{p}$
for $0\leq p\leq\tau-1$, $|u|$ if $p=\tau$, and $q_{p-1}$ for $\tau+1\leq p\leq
k+3$}.
\end{eqnarray*}
Then $|u|$ is in the list of $\mathbf{q}'$ and from
(\ref{INTeq5}), one can easily check that
$\mathcal{P}_{k}(\mathbf{m},\mathbf{q})=\mathcal{P}_{k+1}(\mathbf
{m}',\mathbf{q}')$. Let us notice that this concept, though simple,
will simplify many of our future discussions.

We specify the following values for (\ref{GBasseq1}):
%
\begin{eqnarray}
\label{GBasseq2} \kappa&=&k+1,
\nonumber
\\
n_p&=&\frac{m_p}{1+t} \qquad\mbox{if $0\leq p<\tau$}\quad
\mbox{and}\quad
m_p \qquad\mbox{if $\tau\leq p\leq\kappa$},
\\
\rho_p&=&q_p \qquad\mbox{for $0\leq p\leq\kappa+1$}.
\nonumber
\end{eqnarray}
Let $\lambda=0$. From Theorem~\ref{GBthm2}, it follows that
\[
p_{N,u}\leq2\log2+Y_0-\sum_{0\leq p\leq k+1}m_p
\bigl(\theta(q_{p+1})-\theta(q_p)\bigr).
\]
Let $y_0^1,y_0^2,\ldots,y_{k+1}^1,y_{k+1}^2$ be jointly Gaussian r.v.'s
defined in Theorem~\ref{GBthm2} and be independent of $h$.
For $j=1,2$, we define $(X_p^j)_{0\leq p\leq k+2}$ in the same way as
$(X_p)_{0\leq p\leq k+2}$ by using $k,\mathbf{m},\mathbf{q}$,
and $(y_p^j)_{0\leq p\leq k+1}$.
Since $y_p^1$ and $y_p^2$ are independent of each other for each $\tau
\leq p\leq k+1$, it implies $Y_\tau=X_\tau^1+X_\tau^2$.
To bound $Y_0$ from above, we need the following lemma, which can be
proven by following the same idea as Proposition~\ref{ImportantThm1}
in Section~\ref{sec6} below and is left to the reader.

\begin{lemma}\label{GBlem1}
Suppose that $\eta$ is a constant which takes value $1$ or $-1$.
Consider two jointly Gaussian r.v.'s $y_1$ and $y_2$ such that
$Ey_1^2=Ey_2^2$ and $Ey_1y_2=\eta tEy_1^2$.
Consider two functions $F_1$ and $F_2$ such that their first four
derivatives are uniformly bounded. Then for any values of
$x_1,x_2$ and $m>0$ we have
%
\begin{eqnarray}
\label{GBlem1eq1}
&&
\frac{1+t}{m}\log E\exp\frac{m}{1+t}
\bigl(F_1(x_1+y_1)+F_2(x_2+y_2)
\bigr)
\nonumber\\[-8pt]\\[-8pt]
&&\qquad\leq\sum_{j=1,2}\frac{1}{m}\log E\exp
mF_j(x_j+y_j).
\nonumber
\end{eqnarray}
\end{lemma}

Since $y_p^1$ and $y_p^2$ satisfy $Ey_p^1y_p^2=\eta tE(y_p^1)^2=\eta
tE(y_p^2)^2$ for $0\leq p<\tau$,
using (\ref{GBlem1eq1}) and decreasing induction, $Y_0\leq
X_0^1+X_0^2=2X_0$.
Hence, we conclude that for any given numbers $k,\mathbf{m},\mathbf
{q}$, we can find parameters (\ref{GBasseq1}) such that
$\mathcal{P}_k(\mathbf{m},\mathbf{q})$ can be recovered by the
right-hand side of (\ref{GBthm2eq1}), that is, $p_{N,u}\leq
2\mathcal{P}_k(\mathbf{m},\mathbf{q})$. Now, to prove Proposition
\ref{INTprop1}, we have to find suitable parameters
(\ref{GBasseq1}) for Guerra's bound. It turns out that this can be done
and leads to the following three crucial propositions. First, we have
the following result.

\begin{proposition}\label{MEprop1}
For $0<t\leq1$, there exists a number $\varepsilon^*<0$ depending
only on $t$, $\xi$ and $h$ such that for every $u\leq0$,
$p_{N,u}\leq2\mathcal{P}(\xi,h)-\varepsilon^*$.
\end{proposition}

This proposition means that the overlap takes essentially nonnegative
values, which is mainly due to the presence of the external field, that
is, $Eh^2\neq0$.
Let $\mu$ be the Parisi measure and $c$ be the smallest value of its
support. Recall the definition of $\varphi_v(u,t)$ corresponding to
$\mu$ from (\ref{INTeq4}). Two crucial facts about $Y_0$ that will
be derived in Sections~\ref{sec3} and~\ref{sec5} below are that for
arbitrary choice of (\ref{GBasseq1}), the second partial
derivative of $Y_0$ with respect to $\lambda$ is bounded by $1$ and if
we choose (\ref{GBasseq1}) properly, the first partial derivative
of $Y_0$ at $\lambda=0$ roughly gives the formulation of~$\varphi_c$.
From Guerra's bound and these facts, they imply our next proposition.

\begin{proposition}\label{ME0cprop2}
For $0\leq u\leq c$ and $0\leq t\leq1$ we have
%
\begin{equation}
\label{ME0cthm1eq1} p_{N,u}\leq2\mathcal{P}(\xi,h)-\tfrac{1}{2}
\varphi_c(u,t)^2.
\end{equation}
If $0\leq t<1$, then there exists a $\gamma>0$ depending on the Parisi
measure $\mu$ and $t$ such that
%
\begin{equation}
\label{ME0cthm1eq2} p_{N,u}\leq2\mathcal{P}(\xi,h)-\tfrac{1}{16}
\varphi_c(c,t)^2
\end{equation}
for every $c\leq u\leq c+\gamma$.
\end{proposition}

At last, we investigate the upper bound for $p_{N,u}$ when $u>c'$ for
some fixed \mbox{$c'>c$}. This strongly relies on the assumption that these
two SK models use different disorders, that is, $0<t<1$. Our main
result is stated as follows.

\begin{proposition}\label{MEposprop1}
Suppose that $0<t<1$ and $c<c'<1$. Then there exists $\varepsilon^*>0$
such that
$p_{N,u}\leq2\mathcal{P}(\xi,h)-\varepsilon^*$ for every $c'\leq
u\leq1$, where $\varepsilon^*$ depends only on $t,\xi,h,c'$.
\end{proposition}

These propositions are the main ingredients of the proof of Proposition
\ref{INTprop1} and their proofs are deferred to Sections \ref
{sec4},~\ref{sec5} and~\ref{sec6}, respectively. Now, let us proceed
to prove Proposition~\ref{INTprop1}.

\begin{pf*}{Proof of Proposition~\ref{INTprop1}}
Let $0<t<1$ be fixed. From Proposition~\ref{MEprop1}, there exists
$\varepsilon_1^*$ depending only on $t,\xi$ and $h$
such that for every $-1\leq u\leq0$,
%
\begin{equation}
\label{MEprop1proofeq0} p_{N,u}\leq2\mathcal{P}(\xi,h)-
\varepsilon_1^*.
\end{equation}
Now, for given $\varepsilon>0$, we set
\[
\varepsilon_2^*=\tfrac{1}{2}\min\bigl\{\varphi_c(w,t)^2\dvtx 0
\leq w\leq c, |w-u_t|\geq\varepsilon\bigr\}.
\]
Since $u_t$ is the unique solution of $\varphi_c(\cdot,t)$ in $
[0,c ]$,
it follows that $\varepsilon_2^*>0$ and from (\ref{ME0cthm1eq1}),
%
\begin{equation}
\label{MEprop1proofeq1} p_{N,u}\leq2\mathcal{P}(\xi,h)-
\varepsilon_2^*,
\end{equation}
whenever $0\leq u\leq c$ and $|u-u_t|\geq\varepsilon$.
Since we also know $\varphi_c(c,t)<0$, from (\ref{ME0cthm1eq2}),
there exists some $\gamma>0$ depending only on $\mu$ and $t$ such that
%
\begin{equation}
\label{MEprop1proofeq2} p_{N,u}\leq2\mathcal{P}(\xi,h)-
\varepsilon_3^*
\end{equation}
for every $c\leq u\leq c+\gamma$, where $\varepsilon_3^*=\varphi_c(c,t)^2/16>0$.
Let us put $c'=c+\gamma$ in Proposition~\ref{MEposprop1}. Then
there exists $\varepsilon_4^*>0$ depending only on $t,\xi,h,c'$ such that
%
\begin{equation}
\label{MEprop1proofeq3} p_{N,u}\leq2\mathcal{P}(\xi,h)-
\varepsilon_4^*,
\end{equation}
whenever $c'\leq u\leq1$. Finally, we obtain (\ref{INTprop1eq1})
by combining (\ref{MEprop1proofeq0}),
(\ref{MEprop1proofeq1}), (\ref{MEprop1proofeq2})
and (\ref{MEprop1proofeq3}) together and letting $\varepsilon^*=\min
(\varepsilon_1^*,\varepsilon_2^*,\varepsilon_3^*,\varepsilon_4^* )$.
\end{pf*}


\section{Preliminary results}\label{sec3}

Let $k,\mathbf{m},\mathbf{q}$ be given by (\ref{INTeq2}). Suppose
that\break $(z_p)_{0\leq p\leq k+1}$ are independent Gaussian r.v.'s
with $Ez_p^2=\xi'(q_{p+1})-\xi'(q_p)$. Starting with $A_{k+2}(x)=\log
\C x$, we define
%
\begin{equation}
\label{PReq1} A_p(x)=\frac{1}{m_p}\log E\exp
m_pA_{p+1}(x+z_p)
\end{equation}
for $0\leq p\leq k+1$. If $m_p=0$, we define $A_p(x)=EA_{p+1}(x+z_p)$.
Recall $X_0$ from (\ref{INTeq5}).
It should be clear that
\[
X_p=A_p \biggl(h+\sum_{0\leq n<p}z_n
\biggr)
\]
for every $1\leq p\leq k+2$ and $X_0=EA_0(h)$.
Since we will be working with $(A_p)_{0\leq p\leq k+2}$ for much of the
remainder of this paper, we summarize some quantitative results in
Lemma~\ref{PRlem1}.

\begin{lemma}\label{PRlem1}
For every $0\leq p\leq k+2$, we have
%
\begin{eqnarray}
\label{PRlem1eq1}
A_p(x)&=&A_p(-x),\qquad
\bigl|A_p'\bigr|\leq1,\nonumber\\
\frac{1}{C\C^2 x}&\leq&
A_p''(x)\leq\min\biggl(1,
\frac{C}{\C^2 x} \biggr),
\\
\bigl|A_p^{(3)}\bigr|&\leq&4,\qquad \bigl|A_{p}^{(4)}\bigr|\leq8,
\nonumber
\end{eqnarray}
where $C$ is a constant depending only on $\xi$.
\end{lemma}

\begin{pf} The Poisson--Dirichlet cascade was of great importance in
the study of the random energy model and generalized random energy
model in~\cite{D81,R87} and was put forward to the SK model,
in particular, in~\cite{ASS03,BS98}. Following similar ideas
in these works, it is known from Theorem 14.2.1~\cite{Talag11} that
$A_p$ has a very beautiful representation via the Poisson--Dirichlet
cascade [see (\ref{addeq1}) below]. Our argument will be started
with such representation and is concentrated on the inequality
$A_p''(x)\leq C(\C^2x)^{-1}$. For the other statements, one may refer
to Lemma 14.7.16~\cite{Talag11}. Since $A_p$ is an even function,
it suffices\vadjust{\goodbreak} to prove that $A_p''(x)\leq C\exp(-2x)$ for all $x\geq0$.
Let $\tau_1\geq1$ be the smallest integer with $m_{\tau_1}>0$ and
$\tau_2\leq k$ be the largest integer with $m_{\tau_2}<1$.
Suppose for the moment that there exists $C_1>0$ such that $
A_p''(x)\leq C_1\exp(-2x )$ for all $x\geq0$ and $\tau_1\leq p\leq\tau
_2$. By definition of $A_p$, for all $x\geq0$ and
$0\leq p<\tau_1$, we have that
\[
A_p(x)=EA_{\tau_1} \biggl(x+\sum
_{p\leq n<\tau_1}z_n \biggr)
\]
and then
\begin{eqnarray*}
A_p''(x)&=&EA_{\tau_1}''
\biggl(x+\sum_{p\leq n<\tau_1}z_n \biggr)
\\
&\leq&2C_1\exp(-2x)E\C\biggl(2\sum
_{p\leq n<\tau_1}z_n \biggr)
\\
&=&2C_1\exp\bigl(2 \bigl(\xi'(q_{\tau_1})-
\xi'(q_p) \bigr) \bigr)\exp(-2x).
\end{eqnarray*}
Also for all $x\geq0$ and $\tau_2<p\leq k+2$, it is easy to see that
\[
A_p(x)=\log E\C\biggl(x+\sum_{p\leq n<k+2}z_n
\biggr)=\log\C x+\frac
{1}{2}\bigl(\xi'(1)-
\xi'(q_p)\bigr)
\]
and so $A_p''(x)=(\C^2 x)^{-1}\leq C_2\exp(-2x)$ for some constant $C_2>0$.
If we set
\[
C=\max\bigl(C_2,2C_1\exp\bigl(2\xi'(1)
\bigr) \bigr),
\]
then $A_p''(x)\leq C\exp(-2x)$ for all $x\geq0$ and $0\leq p\leq k+2$.
So in the following, we may assume, without loss of generality, that
$0<m_1,m_{k}<1$ and $1\leq p\leq k$. Also, from the discussion right
below Theorem~\ref{GBthm2}, we may let $0<m_1<m_2<\cdots<m_k<1$.

For $p'$ with $p\leq p'\leq k$ and $j_p,j_{p+1},\ldots,j_{p'-1}\in
\Bbb{N}$, we consider a nonincreasing rearrangement
$(u_{j_pj_{p+1}\cdots j_{p'-1}j})_{j\in\Bbb{N}}$ of a Poisson point
process of intensity measure $x^{-m_{p'}-1}\,dx$.
All of these are independent of each other.
For $\alpha=(j_p,j_{p+1},\ldots,j_{k})\in\Bbb{N}^{k+1-p}$,
we set
\[
u_\alpha^*=u_{j_p}u_{j_pj_{p+1}}\cdots u_{j_p j_{p+1}\cdots j_{k}}
\]
and
\[
\nu_\alpha=\frac{u_\alpha^*}{\sum_\gamma u_\gamma^*}.
\]
This family of random weights is called the Poisson--Dirichlet cascade
associated with the sequence $0<m_p<m_{p+1}< \cdots< m_k<1$.
For each $p'$ with $p\leq p'\leq k$, let us consider a sequence of
independent copies of $z_{p'}$,
\[
(z_{p',j_{p},j_{p+1},\ldots,j_{p'}})_{
j_{p},j_{p+1},\ldots,j_{p'}\in\Bbb{N}}.
\]
These sequences are independent of each other and of
$(u_{j_pj_{p+1}\cdots j_{p'-1}j})_{j\in\Bbb{N}}$
for $p\leq p'\leq k$ and $j_p,j_{p+1},\ldots,j_{p'-1}\in\Bbb{N}$. To
simplify the notation,\vadjust{\goodbreak}
for $\alpha=(j_p,j_{p+1},\ldots,j_{k})\in\Bbb{N}^{k+1-p}$, we write
\[
z_{p',\alpha}=z_{p',j_{p},j_{p+1},\ldots,j_{p'}}.
\]
Then from Theorem 14.2.1~\cite{Talag11},
%
\begin{equation}
\label{addeq1} A_p(x)=E\log\sum_{\alpha}
\nu_\alpha\C(x+z_\alpha)+\frac{1}{2}\bigl(
\xi'(1)-\xi'(q_{k+1})\bigr),
\end{equation}
where
\[
z_\alpha=\sum_{p\leq p'\leq k}z_{p',\alpha}.
\]
Taking derivatives, we obtain
\begin{eqnarray*}
A_p''(x)&=&1-E \biggl(\frac{\sum_\alpha\nu_\alpha\sh(x+z_\alpha
)}{\sum_\alpha\nu_\alpha\C(x+z_\alpha)}
\biggr)^2
\\
&=&E\biggl(\frac{\sum_\alpha\nu_\alpha(\C(x+z_\alpha)-\sh
(x+z_\alpha))}{\sum_\alpha\nu_\alpha\C(x+z_\alpha)}
\\
&&\hspace*{13pt}{} \times\frac{\sum_\alpha\nu_\alpha(\C(x+z_\alpha)+\sh(x+z_\alpha
))}{\sum_\alpha\nu_\alpha\C(x+z_\alpha)}\biggr)
\\
&\leq&2E \biggl(\frac{\sum_\alpha\nu_\alpha\exp(-z_\alpha)}{\sum_\alpha
\nu_\alpha\C(x+z_\alpha)} \biggr)\exp(-x)
\\
&\leq&4E \biggl(\frac{\sum_\alpha\nu_\alpha\exp(-z_\alpha)}{\sum_\alpha
\nu_\alpha\exp(z_\alpha)} \biggr)\exp(-2x),
\end{eqnarray*}
where the first inequality holds since $\C y-\sh y=\exp(-y)$ and
$|{\sh y}|\leq\C y$, while the second
inequality follows from $2\C y\geq\exp y$. Let us now turn to the
computation of this quantity
\[
\gamma_p:=E \biggl(\frac{\sum_\alpha\nu_\alpha\exp(-z_\alpha
)}{\sum_\alpha\nu_\alpha\exp(z_\alpha)} \biggr).
\]
Set $F(x_p,x_{p+1},\ldots,x_{k})=\sum_{p\leq p'\leq k}x_{p'}$ for
$(x_p,x_{p+1},\ldots,x_k)\in\Bbb{R}^{k+1-p}$.
For $\alpha\in\Bbb{N}^{k+1-p}$, define the random variables
\begin{eqnarray*}
F(\alpha)&=&F(z_{p,\alpha},\ldots,z_{k,\alpha}),
\\
U(\alpha)&=&\exp\bigl(-2F(\alpha)\bigr).
\end{eqnarray*}
Then we can write
%
\begin{equation}
\label{PRlem1proofeq2} \gamma_p=E\frac{\sum_\alpha\nu_\alpha U(\alpha
)\exp F(\alpha
)}{\sum_{\alpha}\nu_\alpha\exp F(\alpha)}.
\end{equation}
Starting from
\[
F_{k+1}=F(z_p,z_{p+1},\ldots,z_{k}),\vadjust{\goodbreak}
\]
we define by decreasing induction for $p\leq p'\leq k$,
\[
F_{p'}=\frac{1}{m_{p'}}\log E_{p'}\exp
m_{p'}F_{p'+1},
\]
where $E_{p'}$ means the expectation with respect to the r.v.'s
$z_{p'},z_{p'+1},\ldots,z_{k}$.
We also define for $p\leq p'\leq k$,
\[
W_{p'}=\exp m_{p'}(F_{p'+1}-F_{p'}).
\]
From formula $(14.27)$ in~\cite{Talag11}, (\ref{PRlem1proofeq2})
can be computed as
%
\begin{equation}
\label{PRlem1proofeq1} \gamma_p=EW_pW_{p+1}
\cdots W_{k}\exp(-2F_{k+1} ).
\end{equation}
Using the independence of $z_p,z_{p+1},\ldots,z_{k}$, it is easy to
compute that
\[
F_{p'}=\sum_{n=p}^{p'-1}z_n+
\frac{1}{2}\sum_{n=p'}^{k}m_{n}
\bigl(\xi'(q_{n+1})-\xi'(q_{n})
\bigr).
\]
Therefore, we obtain
\[
W_{p'}=\exp\biggl(m_{p'}z_p'-
\frac{m_{p'}^2}{2}\bigl(\xi'(q_{p'+1})-\xi'(q_{p'})
\bigr) \biggr)
\]
and from (\ref{PRlem1proofeq1}), this implies
\begin{eqnarray*}
\gamma_p&=&E\exp\Biggl(\sum_{p'=p}^{k}(m_{p'}-2)z_{p'}-
\frac
{1}{2}\sum_{p'=p}^{k}m_{p'}^2
\bigl(\xi'(q_{p'+1})-\xi'(q_{p'})
\bigr) \Biggr)
\\
&=&\exp\Biggl(\frac{1}{2}\sum_{p'=p}^{k}
\bigl((m_{p'}-2)^2-m_{p'}^2 \bigr)
\bigl(\xi'(q_{p'+1})-\xi'(q_{p'})
\bigr) \Biggr)
\\
&=&\exp\Biggl(2\sum_{p'=p}^{k}
(1-m_{p'} ) \bigl(\xi'(q_{p'+1})-
\xi'(q_{p'})\bigr) \Biggr)
\\
&\leq&\exp\bigl(2\xi'(1)\bigr).
\end{eqnarray*}
Finally, we are done by letting $C=4\exp(2\xi'(1))$.
\end{pf}

As a consequence of Lemma~\ref{PRlem1}, we have the following lemma.

\begin{lemma}\label{PRlem3}
There exists a number $M$ depending only on $\xi$ and $h$ such that
for every $0\leq p\leq k+2$,
%
\begin{eqnarray}
\label{PRlem3eq1} EA_p'\bigl(h+
\chi_1'\bigr)A_p'\bigl(h+
\chi_2'\bigr)&\geq&\frac{1}{M},
\\
\label{PRlem3eq2} EA_p''\bigl(h+
\chi_1''\bigr)A_p''
\bigl(h+\chi_2''\bigr)&\geq&
\frac{1}{M},
\end{eqnarray}
where $\chi_1'$, $\chi_2',\chi_1'',\chi_2''$ are jointly Gaussian
r.v.'s with the same variance $\xi'(q_p)$ and $E\chi_1'\chi_2'=0$
independent of $h$.
\end{lemma}

\begin{pf}
The first inequality is Lemma 14.12.8~\cite{Talag11} and from there a
similar argument yields the second inequality.
\end{pf}

Recall that the external field $h$ in this paper is always assumed to
satisfy \mbox{$Eh^2>0$}.
Based on this assumption, we set up the definition for the Parisi measure.

\begin{definition}\label{PRdef1}
Given $\varepsilon>0$, we say that $k,\mathbf{m},\mathbf{q}$ satisfy
condition $\operatorname{MIN}(\varepsilon)$
if the following occurs. First, the sequences
\begin{eqnarray*}
\mathbf{m}&=&(m_0,m_1,\ldots,m_k,m_{k+1}),
\\
\mathbf{q}&=&(q_0,q_1,\ldots,q_{k+1},q_{k+2})
\end{eqnarray*}
satisfy
\begin{eqnarray*}
m_0&=&0<m_1<\cdots<m_{k}<m_{k+1}=1,
\\
q_0&=&0<q_1<\cdots<q_{k+1}<q_{k+2}=1.
\end{eqnarray*}
In addition,
\[
\mathcal{P}_k(\mathbf{m},\mathbf{q})\leq\mathcal{P}(\xi,h)+
\varepsilon
\]
and
\[
\mbox{$\mathcal{P}_k(\mathbf{m},\mathbf{q})$ realizes the minimum
of $\mathcal{P}_k$ over all choices of $\mathbf{m}$ and $
\mathbf{q}$}.
\]
\end{definition}

\begin{definition}\label{sec3def1}
Suppose that $\mu$ is a probability measure associated to $k,\mathbf
{m},\mathbf{q}$. Then we say that $\mu$ is $\varepsilon$-stationary
for some $\varepsilon>0$ if $k,\mathbf{m},\mathbf{q}$ satisfy
condition $\operatorname{MIN}(\varepsilon)$.
\end{definition}

Let us note from Lemma 14.5.5~\cite{Talag11} that for any given
$\varepsilon>0$, we can find an $\varepsilon$-stationary measure $\mu
$ associated to some $k,\mathbf{m},\mathbf{q}$.

\begin{definition}\label{sec3def2}
We say that a probability measure $\mu$ is a Parisi measure
(corresponding to the function $\xi$ and external field $h$)
if there exist a sequence $(\varepsilon_n)$ with $\varepsilon_n
\downarrow0$ and a sequence of probability measures $(\mu_n)$ such
that the following two conditions hold:
\begin{eqnarray*}
&\mbox{$\mu_n$ is $\varepsilon_n$-stationary},&
\\
&\mbox{$\mu$ is the limit of $(\mu_n)$.}& 
\end{eqnarray*}
\end{definition}

Definition~\ref{PRdef1} is the same as Definition 14.5.3
\cite{Talag11}, while our definition of the stationarity in Definition
\ref{sec3def1} is stronger than that in Definition 14.11.4~\cite{Talag11}.
This is for technical purposes and it should be clear
that under these assumptions, our future arguments are still valid.\vadjust{\goodbreak}

\begin{lemma}\label{PRlem4}
Suppose $h\neq0$ and $k,\mathbf{q},\mathbf{m}$ satisfy condition
$\operatorname{MIN}(\varepsilon)$. Then for $1\leq p\leq k+1$,
%
\begin{eqnarray}
\label{PRlem4eq1} E W_1\cdots W_{p-1}A_{p}'(
\zeta_p)^2&=&q_p,
\\
\label{PRlem4eq2} \xi''(q_p)E
W_1\cdots W_{p-1}A_p''(
\zeta_p)^2&\leq&1+M\varepsilon^{1/6},
\end{eqnarray}
where $\zeta_p=h+\sum_{0\leq n<p}z_n$
and
\[
W_p=\exp m_p\bigl(A_{p+1}(
\zeta_{p+1})-A_p(\zeta_p)\bigr)=\exp
m_p(X_{p+1}-X_p).
\]
Here, $M$ is a constant depending only on $\xi$ and $h$.
\end{lemma}

\begin{pf}
These results are $(14.222)$ and $(14.461)$ in~\cite{Talag11}.
\end{pf}

At the end of this section we will find a manageable bound for
$p_{N,u}$ via Guerra's bound.
Recall that the right-hand side of (\ref{GBthm2eq1}) depends on
(\ref{GBasseq1}). If we keep every parameter but $\lambda$ fixed,
then it is a quantity depending only on $\lambda$ and, for clarity, we
denote it by $\alpha(\lambda)$. For the same reason, we also think of $Y_0$
as a function of~$\lambda$.
Recall the r.v.'s $(y_p^j)_{0\leq p\leq\kappa,j=1,2}$ defined in
Theorem~\ref{GBthm2}.
Suppose that $(y_p)_{0\leq p\leq\kappa}$ are independent Gaussian
r.v.'s with $E(y_p)^2=\xi'(\rho_{p+1})-\xi'(\rho_p)$ for $0\leq
p\leq\kappa$.
Starting with
\[
D_{\kappa+1}(x)=\log\C x,
\]
we define $D_p$ for $0\leq p\leq\kappa$ by decreasing induction:
\[
D_p(x)=\cases{ %
\displaystyle \frac{1}{n_p}\log E_p
\exp n_p D_{p+1}(x+y_p),&\quad if $\tau\leq p\leq
\kappa$,
\cr
\displaystyle \frac{1}{(1+t)n_p}\log E_p\exp(1+t)n_pD_{p+1}(x+y_p),&\quad
if $0\leq p<\tau$,}
\]
where $E_p$ means the expectation with respect to $y_{n}$ for $p\leq
n\leq\kappa$. If $n_p=0$ for some $p$,
then we define $D_p(x)=E_{p}D_{p+1}(x+y_p)$.
For $j=1,2$ and $1\leq p\leq\kappa+1$, set
\[
\zeta_p^j=h+\sum_{0\leq n<p}y_n^j.
\]

\begin{proposition}\label{prop3}
If $n_p=0$ for every $0\leq p<\tau$, then
%
\begin{eqnarray}
\label{prop3eq3} Y_0(0)&=&ED_\tau\bigl(
\zeta_\tau^1\bigr)+ED_\tau\bigl(
\zeta_\tau^2\bigr),
\\
\label{prop3eq4} Y_0'(0)&=&ED_\tau'
\bigl(\zeta_\tau^1\bigr)D_\tau'
\bigl(\zeta_\tau^2\bigr).
\end{eqnarray}
For the second derivative of $Y_0$, we have for every $\lambda$,
%
\begin{equation}
\label{prop3eq5} 0\leq Y_0''(\lambda)
\leq1.
\end{equation}
\end{proposition}

\begin{pf}
The proofs of (\ref{prop3eq3}) and (\ref{prop3eq4}) are
essentially the same as that of part~(b) of Proposition 14.6.4 \cite
{Talag11}.
Also, (\ref{prop3eq5}) and Lemma 14.6.5~\cite{Talag11} have the
same proof.\vadjust{\goodbreak}
\end{pf}

\begin{corollary}\label{PRcor1} We have
%
\begin{equation}
\label{PRcor1eq1} p_{N,u}\leq\inf_\lambda\alpha(\lambda)\leq
\alpha(0)-\frac
{1}{2}\alpha'(0)^2.
\end{equation}
\end{corollary}

\begin{pf}
This is an immediate consequence of (\ref{prop3eq5}).
\end{pf}

Let us remark here that (\ref{PRcor1eq1}) helps us in at least two
ways: First, it reduces the difficulty of choosing parameters
since we do not have to choose $\lambda$ now. Second, this inequality
gives us a reasonable way to choose parameters.
Roughly speaking, in many cases, we choose parameters in such a way
that the quantity $\alpha(0)$ is very close to $\mathcal{P}(\xi,h)$, while
the term $\alpha'(0)^2/2$ is the error that we expect to obtain on the
right-hand side of (\ref{INTprop1eq1}).


\section{\texorpdfstring{Proof of Proposition \protect\ref{MEprop1}}{Proof of Proposition 3}}\label{sec4}

This section is devoted to proving Proposition~\ref{MEprop1}.
Our approach is based on Talagrand's proof of the positivity of the
overlap in Section 14.12~\cite{Talag11}.
Suppose that $u=-v$ for $0\leq v\leq1$. Proposition~\ref{MEprop1}
relies on the following two results:

\begin{proposition}\label{MEprop2}
There exists $\delta>0$ and $\varepsilon_0>0$ depending only on $\xi
$ and $h$ with the following
property. Whenever we can find $k,\mathbf{m},\mathbf{q}$ that satisfy
condition $\operatorname{MIN}(\varepsilon_0)$ and for an
integer $s$ with $1\leq s\leq k+1$,
\[
m_{s-1}\leq\delta\quad\mbox{and}\quad q_s\geq v-\delta,
\]
then we can find parameters in (\ref{GBthm2eq1}) such that
$p_{N,u}\leq2\mathcal{P}(\xi,h)-1/M$,
where $M$ depends only on $\xi$ and $h$.
\end{proposition}

\begin{proposition}\label{MEprop3}
Consider $\delta$ as in Proposition~\ref{MEprop2}. Then we can find
$\varepsilon_1>0$ with the following property. Whenever we can find
$k,\mathbf{m},\mathbf{q}$ such that $\mathcal{P}_k(\mathbf
{m},\mathbf{q})\leq\mathcal{P}(\xi,h)+\varepsilon_1$ and an
integer $s$
with $1\leq s\leq k+1$,
\[
m_s\geq\delta\quad\mbox{and}\quad q_s\leq v-\delta,
\]
then we can find parameters in (\ref{GBthm2eq1}) such that
$p_{N,u}\leq2\mathcal{P}(\xi,h)-1/M$,
where $M$ depends only on $\xi$ and $h$.
\end{proposition}

\begin{pf*}{Proof of Proposition~\ref{MEprop1}}
Let $v\geq0$. Consider $\delta$, $\varepsilon_0$ as in Proposition~\ref{MEprop2}
and $\varepsilon_1$ as in Proposition~\ref{MEprop3}.
Suppose that $k,\mathbf{m},\mathbf{q}$ is a triplet satisfying
$\operatorname{MIN}(\min(\varepsilon_0,\varepsilon_1))$. Here, the existence of
such $k,\mathbf{m},\mathbf{q}$ is ensured by Lem\-ma~14.5.5~\cite
{Talag11}. Let $1\leq s\leq k+1$ be the largest integer such that
$m_{s-1}\leq\delta$. If $q_s\geq v-\delta$, we apply Proposition
\ref{MEprop2}. Otherwise we have $q_s\leq v-\delta$. If $s=k+1$,
then $m_s=m_{k+1}=1\geq\delta$.
If $s<k+1$, then from the definition of $s$, $m_s\geq\delta$. In both
cases, we conclude Proposition~\ref{MEprop1}
by using Proposition~\ref{MEprop3} and we are done.
\end{pf*}

Note that since the proof of Proposition~\ref{MEprop3} is
essentially the same as that of Proposition~\ref{MEposprop1},
we defer it to Section~\ref{sec6}. Now we turn to the proof of
Proposition~\ref{MEprop2} and proceed with the following lemma:\vadjust{\goodbreak}

\begin{lemma}\label{PRlem2}
Suppose that $A\dvtx \Bbb{R}\rightarrow\Bbb{R}$ has uniformly bounded
first and second derivatives.
Consider two independent pairs of jointly Gaussian r.v.'s $(\chi_1,\chi
_2)$ and $(\chi_1',\chi_2')$, all of variance $a$,
and a standard Gaussian r.v. $\chi$. These r.v.'s are independent of
$h$. Then we have
%
\begin{eqnarray}
\label{PRlem2eq1}
&&
\bigl\llvert EA'(h+\chi_1)A'(h+
\chi_2)-EA'\bigl(h+\chi_1'
\bigr)A'\bigl(h+\chi_2'\bigr)\bigr\rrvert
\nonumber\\[-8pt]\\[-8pt]
&&\qquad\leq\bigl|E\chi_1\chi_2-E\chi_1'
\chi_2'\bigr|EA''(h+\chi
\sqrt{a})^2.
\nonumber
\end{eqnarray}
\end{lemma}

\begin{pf}
This is a typical application of the Gaussian interpolation technique
and the Cauchy--Schwarz inequality. For details, one may refer to
Lem\-ma~14.9.5~\cite{Talag11}.
\end{pf}

Suppose that $k$, $\mathbf{m}$, $\mathbf{q}$ is a triplet satisfying
$\operatorname{MIN}(\varepsilon)$.
Based on our discussion in Section~\ref{sec1}, we may assume, without
loss of generality, that $v=q_a$ for some $a$.
The only thing we have to keep in mind is that when using
(\ref{PRlem4eq2}), we will not be able to use the value $p=a$.
From the assumption that $q_s\geq v-\delta$, we divide our discussion
into two cases $v-\delta\leq q_s\leq v$ and $q_s> v$.
First, let us proceed with the case that for an integer $s$ with $1\leq
s\leq k+1$,
%
\begin{equation}
\label{eq6} m_{s-1}\leq\delta\quad\mbox{and}\quad v-\delta\leq q_s
\leq v.
\end{equation}
Note that $s\leq a$. We consider the following numbers:
%
\begin{eqnarray}
\label{eq10} \quad\tau&=&1,
\nonumber
\\
\kappa&=&k+2-a,
\nonumber\\[-8pt]\\[-8pt]
n_0&=&0,\qquad n_1=m_a,\qquad n_2=m_{a+1},\ldots,\qquad n_{\kappa}=m_{k+1}=1,
\nonumber
\\
\rho_0&=&0,\qquad\rho_1=v=q_a,\qquad
\rho_2=q_{a+1},\ldots,\qquad\rho_{\kappa+1}=q_{k+2}=1
\nonumber
\end{eqnarray}
and apply (\ref{eq10}) to Theorem~\ref{GBthm2}. Recall that we
use $\alpha$ to denote the right-hand side of (\ref{GBthm2eq1}).

\begin{lemma}\label{lem3} Assuming (\ref{eq6}) and (\ref{eq10}),
we have
%
\begin{equation}
\label{lem3eq1} \alpha(0)\leq2\mathcal{P}_k(\mathbf{m},
\mathbf{q})+M\delta.
\end{equation}
\end{lemma}

\begin{pf}
The proof is essentially the same as that of Lemma 14.12.7 in~\cite{Talag11}.
\end{pf}

In view of (\ref{PRcor1eq1}) and (\ref{lem3eq1}), our goal is
then to bound $\alpha'(0)$ from below. Proposition~\ref{prop3}
implies that $D_1(x)=A_{a}(x)$
and so
%
\begin{equation}
\label{eq2} \alpha'(0)=EA_a'(h+
\chi_1)A_a'(h+\chi_2)+v,
\end{equation}
where $\chi_1$ and $\chi_2$ are Gaussian with $E(\chi_1)^2=E(\chi
_2)^2=\xi'(v)$ and $E\chi_1\chi_2=-t\xi'(v)$ independent of $h$.
Consider two independent Gaussian\vadjust{\goodbreak} r.v.'s
$\chi_1'$ and $\chi_2'$ with $E(\chi_1')^2=E(\chi_2')^2=\xi'(v)$
independent of $h$. By using
(\ref{PRlem2eq1}),
\begin{eqnarray*}
&&
EA_a'(h+\chi_1)A_a'(h+
\chi_2)
\\
&&\qquad\geq EA_a'\bigl(h+\chi_1'
\bigr)A_a'\bigl(h+\chi_2'
\bigr)-t\xi'(v)EA_a''
\bigl(h+\chi\sqrt{\xi'(v)}\bigr)^2,
\end{eqnarray*}
where $\chi$ is standard Gaussian independent of $h$. Since $\xi
'(v)\leq v\xi''(v)$, it follows that from (\ref{eq2}),
%
\begin{eqnarray}
\label{eq3}\quad \alpha'(0) &\geq& EA_a'
\bigl(h+\chi_1'\bigr)A_a'
\bigl(h+\chi_2'\bigr)+v \bigl(1-t\xi''(v)EA_a''
\bigl(h+\chi\sqrt{\xi(v)}\bigr)^2 \bigr)
\nonumber
\\
&=&EA_a'\bigl(h+\chi_1'
\bigr)A_a'\bigl(h+\chi_2'
\bigr)+v(1-t)
\\
&&{} +tv \bigl(1-\xi''(v)EA_a''
\bigl(h+\chi\sqrt{\xi(v)}\bigr)^2 \bigr).
\nonumber
\end{eqnarray}
To use (\ref{eq3}), we have to bound the quantity
\[
\xi''(v)EA_a''
\bigl(h+\chi\sqrt{\xi'(v)}\bigr)
\]
from above. The starting point of the proof is that from (\ref{PRlem4eq2}),
%
\begin{equation}
\label{eq4} \xi''(q_s)EW_1
\cdots W_{s-1}A_s''(
\zeta_s)^2\leq1+M\varepsilon^{1/6},
\end{equation}
where $\zeta_p=h+\sum_{0\leq n<p}z_n$ and
$W_p=\exp m_p(A_{p+1}(\zeta_{p+1})-A_p(\zeta_p))$.

\begin{lemma}\label{lem5}
Assuming (\ref{eq6}), there exists $\delta_0>0$ depending only on
$\xi$ and $h$ such that when $\delta\leq\delta_0$,
we have
%
\begin{equation}
\label{lem5eq1}\quad \xi''(v)EA_a''
\bigl(h+\chi\sqrt{\xi'(v)}\bigr)^2\leq
\xi''(q_s)EW_1\cdots
W_{s-1}A_s''(
\zeta_s)^2+M\sqrt{\delta}.
\end{equation}
\end{lemma}

\begin{pf}
This is Lemma 14.12.9 in~\cite{Talag11}.
\end{pf}

As a conclusion, by assuming (\ref{eq6}) and using (\ref{eq10}),
we see that (\ref{PRlem3eq1}), (\ref{eq3}), (\ref{eq4}) and
(\ref{lem5eq1}) together imply
%
\begin{equation}
\label{eq8} \alpha'(0)\geq\frac{1}{M}-M\varepsilon^{1/6}-M
\sqrt{\delta}
\end{equation}
for $\delta\leq\delta_0$.

Next, let us consider the other case that for some $1\leq s\leq k+1$,
%
\begin{equation}
\label{eq11} m_{s-1}\leq\delta\quad\mbox{and}\quad q_s>v=q_a.
\end{equation}
Since $q_{a+1}\geq q_a\geq v-\delta$ and $m_{a}\leq m_{s-1}\leq\delta$,
we may assume, without loss of generality, that $s=a+1$. Consider the
following numbers:
%
\begin{eqnarray}
\label{eq9}\qquad \tau&=&1,
\nonumber
\\
\kappa&=&k+2-a,
\nonumber\\[-8pt]\\[-8pt]
n_0&=&0,\qquad n_1=0,\qquad n_2=m_{a+1},\ldots,\qquad
n_\kappa=m_{k+1}=1,
\nonumber
\\
\rho_0&=&0,\qquad \rho_1=v=q_a,\qquad
\rho_2=q_{a+1},\ldots,\qquad\rho_{\kappa+1}=q_{k+2}=1
\nonumber
\end{eqnarray}
and apply (\ref{eq9}) to (\ref{GBthm2eq1}).

\begin{lemma}\label{lem6}
Assuming (\ref{eq11}) and (\ref{eq9}), we have
%
\begin{equation}
\label{lem6eq1} \alpha(0)\leq2\mathcal{P}_k(\mathbf{m},
\mathbf{q})+M\delta.
\end{equation}
\end{lemma}

\begin{pf}
A similar proof as Lemma~\ref{lem3} yields the announced statement.
\end{pf}

Again, our goal is to bound $\alpha'(0)$ from below. From
(\ref{prop3eq4}), we have $D_2(x)=A_{a+1}(x)$ and then
%
\begin{equation}
\label{MEeq1} \alpha'(0)=EA_{a+1}'(h+
\chi_1)A_{a+1}'(h+\chi_2)+v,
\end{equation}
where $\chi_1$ and $\chi_2$ are jointly Gaussian with $E(\chi
_1)^2=E(\chi_2)^2=\xi'(q_{a+1})$ and
$E\chi_1\chi_2=-t\xi'(v)$ independent of $h$. Let $\chi_1'$ and
$\chi_2'$ be two independent Gaussian r.v.'s with $E(\chi_1')^2=E(\chi
_2')^2=\xi'(q_{a+1})$ independent of $h$.
Using (\ref{PRlem2eq1}), we obtain
%
\begin{eqnarray}
\label{MEeq2}
&&
EA_{a+1}'(h+\chi_1)A_{a+1}'(h+
\chi_2)
\nonumber\\
&&\qquad\geq EA_{a+1}'\bigl(h+\chi_1'
\bigr)A_{a+1}'\bigl(h+\chi_2'
\bigr)\\
&&\qquad\quad{}-t\xi'(v)EA_{a+1}''
\bigl(h+\chi\sqrt{\xi'(q_{a+1})}\bigr)^2,
\nonumber
\end{eqnarray}
where $\chi$ is standard Gaussian independent of $h$. Let us apply
$p=a+1$ to (\ref{PRlem4eq2}) and use the fact
$q_{a+1}\geq v$. Then we have
%
\begin{eqnarray}
\label{MEeq3}\quad \xi''(v)EW_1\cdots
W_aA_{a+1}''(
\zeta_{a+1})^2 &\leq& \xi''(q_{a+1})EW_1
\cdots W_aA_{a+1}''(
\zeta_{a+1})^2
\nonumber\\[-8pt]\\[-8pt]
&\leq& 1+M\varepsilon^{1/6}.
\nonumber
\end{eqnarray}

\begin{lemma}\label{MElem1}
Assuming (\ref{eq11}), we have
%
\begin{equation}
\label{MElem1eq1} E|W_1\cdots W_{s-1}-1|\leq M\delta.
\end{equation}
\end{lemma}
\begin{pf}
One can find the proof from Lemma 14.12.9~\cite{Talag11}.
\end{pf}

Using (\ref{MElem1eq1}) and $EA_{a+1}''(\zeta
_{a+1})^2=EA_{a+1}''(h+\chi\sqrt{\xi'(q_{a+1})})^2$,
it follows that from (\ref{MEeq3}),
%
\begin{equation}
\label{MEeq5} \xi''(v)EA_{a+1}''
\bigl(h+\chi\sqrt{\xi'(q_{a+1})}\bigr)^2
\leq1+M\delta+M\varepsilon^{1/6}
\end{equation}
and from (\ref{PRlem3eq1}), (\ref{MEeq1}), (\ref{MEeq2}),
(\ref{MEeq5}) and $\xi'(v)\leq v\xi''(v)$, we then have
%
\begin{eqnarray}
\label{MEeq4} \alpha'(0)&\geq& EA_{a+1}'
\bigl(h+\chi_1'\bigr)A_{a+1}'
\bigl(h+\chi_2'\bigr)+v
\nonumber
\\
&&{} -t\xi'(v)EA_{a+1}''
\bigl(h+\chi\sqrt{\xi'(q_{a+1})}\bigr)^2
\nonumber\\[-8pt]\\[-8pt]
&\geq& \frac{1}{M}+v(1-t)+tv \bigl(1-\xi''(v)EA_{a+1}''
\bigl(h+\chi\sqrt{\xi'(q_{a+1})}\bigr)^2
\bigr)
\nonumber
\\
&\geq& \frac{1}{M}-M\delta-M\varepsilon^{1/6}.
\nonumber
\end{eqnarray}

\begin{pf*}{Proof of Proposition~\ref{MEprop2}}
First we complete the proof for the case (\ref{eq6}). Let $M_1$ be
the constant obtained from (\ref{lem3eq1}) and
(\ref{eq8}) and assume, without loss of generality, that $M_1\geq1$
and $1/16M_1^4\leq\delta_0$. Set $\delta=1/16M_1^4$.
If $\varepsilon\leq\varepsilon_1=(1/4M_1^2)^6$, (\ref{eq8}) implies
\[
\alpha'(0)\geq\frac{1}{M_1}-\frac{M_1}{4M_1^2}-
\frac
{M_1}{4M_1^2}=\frac{1}{2M_1}
\]
and combining this with (\ref{lem3eq1}) yields
\[
\inf_\lambda\alpha(\lambda)\leq2\mathcal{P}_k(\mathbf{m},
\mathbf{q})+M_1\delta-\frac{1}{8M_1^2}\leq2\mathcal{P}(\xi,h)+2
\varepsilon_0 -\frac{1}{16M_1^2}.
\]
Letting $\varepsilon_0$ be sufficiently small completes our proof of
this case.
For the second case (\ref{eq11}), using (\ref{MEeq4}) and Lemma
\ref{lem6}, we may argue similarly to obtain
the announced result.
\end{pf*}


\section{\texorpdfstring{Proofs of Propositions \protect\ref{ME0cprop1} and \protect\ref{ME0cprop2}}
{Proofs of Propositions 1 and 4}}\label{sec5}

Given $0\leq v<1$, recall the definition of $\varphi_v$ from
(\ref{INTeq4}).
In this section we first study how the Guerra bound relates to $\varphi
_v$ and then study some of its basic properties to
conclude Propositions~\ref{ME0cprop1} and~\ref{ME0cprop2}.

Let $k,\mathbf{m},\mathbf{q}$ be given by (\ref{INTeq2}).
Suppose that $\mu$ is the probability measure associated to $k,\mathbf
{m},\mathbf{q}$
and $\Phi$ is the corresponding solution of (\ref{INTeq3}). Recall
the definition of $(A_p)_{0\leq p\leq k+2}$ from (\ref{PReq1}).
Then $\Phi$ and $(A_p)_{0\leq p\leq k+2}$ can be related in the
following way.
Let $(g_p)_{0\leq p\leq k+1}$ be i.i.d. standard Gaussian r.v.'s.
For $q\in[0,1 ]$, we have that $\Phi(x,1)=A_{k+2}(x)$ if
$q=1$ and
\[
\Phi(x,q)=\frac{1}{m_p}\log E\exp m_p A_{p+1}
\bigl(x+g_p\sqrt{\xi'(q_{p+1})-
\xi'(q)} \bigr),
\]
if $q_p\leq q<q_{p+1}$ for some $0\leq p\leq k+1$. In particular, for
$0\leq p\leq k+2$,
%
\begin{equation}
\label{PRadd} \Phi(x,q_p)=A_p(x).
\end{equation}
For fixed $u$ and $v$ with $0\leq u\leq v<1$, we suppose $q_a\leq
v<q_{a+1}$ for some $0\leq a\leq k+1$ and consider numbers
%
\begin{eqnarray}
\label{ME0cass2} \tau&=&1,
\nonumber
\\
\kappa&=&k+3-a,
\nonumber\\
n_0&=&0,\qquad n_1=0,\qquad n_2=m_a,\nonumber\\[-8pt]\\[-8pt]
n_3&=&m_{a+1},
\ldots,\qquad n_\kappa=m_{k+1}=1,
\nonumber
\\
\rho_0&=&0,\qquad\rho_1=u,\qquad \rho_2=v,\nonumber\\
\rho_3&=&q_{a+1},\ldots,\qquad\rho_{\kappa
+1}=q_{k+2}=1.
\nonumber
\end{eqnarray}
Let us apply (\ref{ME0cass2}) to (\ref{GBasseq1}) and recall
that we use $\alpha(\lambda)$ to denote the right-hand side of
(\ref{GBthm2eq1}). Recall that $c$ is the smallest value of the support
of the Parisi measure. Since $Eh^2\neq0$, the positivity of the
overlap implies $c>0$.

\begin{lemma}\label{ME0clem1} For $0<\delta<c$, we have
%
\begin{equation}
\label{ME0clem1eq3} \alpha(0)\leq2\mathcal{P}_k(\mathbf{m},
\mathbf{q})+\mu\bigl( [0,c-\delta]\bigr)\theta(1)+ \bigl(\theta
(v)-\theta(c-
\delta) \bigr)_+.
\end{equation}
The derivative of $\alpha$ at $0$ can be computed as
%
\begin{equation}
\label{ME0clem1eq2} \alpha'(0)=E\frac{\partial\Phi}{\partial x}(h+
\chi_1,v)\,\frac
{\partial\Phi}{\partial x}(h+\chi_2,v)-u,
\end{equation}
where $\chi_1$ and $\chi_2$ are two Gaussian r.v.'s with $E(\chi
_1)^2=E(\chi_2)^2=\xi'(v)$ and $E\chi_1\chi_2=t\xi'(u)$
independent of $h$.
\end{lemma}

\begin{pf}
Without loss of generality, we may assume that $v=q_a$ and $u=q_b$ with
$0\leq b\leq a$.
Let us write
%
\begin{eqnarray}
\label{ME0clem1proofeq1}\quad \sum_{1\leq p\leq\kappa}n_p
\bigl(\theta(\rho_{p+1})-\theta(\rho_p)\bigr)&=&\sum
_{a\leq p\leq k+1}m_p\bigl(\theta(q_{p+1})-
\theta(q_p)\bigr)
\nonumber\\[-8pt]\\[-8pt]
&=&\sum_{1\leq p\leq k+1}m_{p}\bigl(
\theta(q_{p+1})-\theta(q_p)\bigr)-C,
\nonumber
\end{eqnarray}
where
\[
C=\sum_{1\leq p\leq a-1}m_{p} \bigl(
\theta(q_{p+1})-\theta(q_p) \bigr).
\]
If $q_a\leq c-\delta$, then
\begin{eqnarray*}
C&\leq&\max\{m_p\dvtx q_p\leq c-\delta\}\sum
_{0\leq p\leq
a-1}\bigl(\theta(q_{p+1})-\theta(q_p)
\bigr)
\\
&\leq&\mu\bigl( [0,c-\delta] \bigr)\theta(1);
\end{eqnarray*}
if $q_a>c-\delta$, then
\begin{eqnarray*}
C&\leq&\max\{m_p\dvtx q_{p}\leq c-\delta\}\sum
_{0\leq p\leq
a-1}\bigl(\theta(q_{p+1})-\theta(q_p)
\bigr)
\\
&&{} +\sum_{0\leq p\leq a-1:q_{p}>c-\delta}\theta(q_{p+1})-\theta
(q_p)
\\
&\leq&\mu\bigl( [0,c-\delta] \bigr)\theta(1)+\theta(v)-\theta(c-\delta).
\end{eqnarray*}
So (\ref{ME0clem1eq3}) holds.
From (\ref{prop3eq3}), we have $D_2(x)=A_a(x)$ and, consequently,
$Y_0=2EA_{a} (h+\chi)$,
where $\chi$ is Gaussian with $E\chi^2=\xi'(q_a)$. Since $\chi$ has
the same distribution as $\sum_{0\leq p<a}z_p$,
from Jensen's inequality, $A_p(x)\geq EA_{p+1}(x+z_p)$ and iterating
this inequality implies
\[
EA_{a} \biggl(h+\sum_{0\leq p<a}z_p
\biggr)\leq EA_0(h).
\]
So $Y_0\leq2EA_0(h)=2X_0$ and this together with
(\ref{ME0clem1proofeq1}) yields (\ref{ME0clem1eq3}).
Next, using (\ref{prop3eq4}) and (\ref{PRadd}), we obtain
\begin{eqnarray*}
Y_0'(0)&=&EA_a'(h+
\chi_1)A_a'(h+\chi_2)
\\
&=&E\frac{\partial\Phi}{\partial x}(h+\chi_1,q_a)\,
\frac{\partial
\Phi}{\partial x}(h+\chi_2,q_a),
\end{eqnarray*}
where $\chi_1$ and $\chi_2$ are jointly Gaussian with $E(\chi
_1)^2=E(\chi_2)^2=\xi'(q_a)$ and $E\chi_1\chi_2=t\xi'(q_b)$
independent of $h$. This completes our proof.
\end{pf}

Now, suppose that $\mu$ is a Parisi measure and $c$ is the smallest
value of its support. By Definition~\ref{sec3def2}, $\mu$ is the limit
of a sequence of $\varepsilon_n$-stationary measures $(\mu_n)$ such
that $\mathcal{P}(\xi,h,\mu_n)\rightarrow\mathcal{P}(\xi,h)$. By
Definition~\ref{sec3def1}, for each $\mu_n$, there exist
$k,\mathbf{m},\mathbf{q}$ satisfying
$\operatorname{MIN}(\varepsilon_n)$. Here, to clarify notation, we keep
the dependence of $k,\mathbf {m},\mathbf{q}$, and $\varepsilon_n$ on
$n$ implicit. For $u$ and $v$ satisfying $0\leq u\leq v<1$, we consider
numbers (\ref{ME0cass2}) associated to $u,v$ and $\mu_n$, and we use
$\alpha_n$ to denote the right-hand side of~(\ref{GBthm2eq1}). Suppose
that $\Phi_n$ is the solution of (\ref{INTeq3}) associated to $\mu_n$.
Recall that we define $\Phi$ as the uniform limit of $(\Phi_n)$. An
argument similar to the proof of Theorem 3.2~\cite{Talag06} implies
that in the sense of uniform convergence,
\[
\frac{\partial^i\Phi}{\partial x^i}=\lim_{n\rightarrow\infty
}\frac{\partial^i\Phi_n}{\partial x^i}
\]
on $\Bbb{R}\times[0,1 ]$ for $i=1,2,3$.

\begin{proposition}\label{ME0cprop3}
For any $u$ and $v$ satisfying $0\leq u\leq v<1$, we have
%
\begin{equation}
\label{ME0cprop3eq1} \limsup_{n\rightarrow\infty}\alpha_n(0)\leq2
\mathcal{P}(\xi,h)+\bigl(\theta(v)-\theta(c)\bigr)_+
\end{equation}
and
%
\begin{equation}
\label{ME0cprop3eq2} \lim_{n\rightarrow\infty}\alpha_n'(0)=E
\frac{\partial\Phi
}{\partial x}(h+\chi_1,v)\,\frac{\partial\Phi}{\partial x}(h+
\chi_2,v)-u,
\end{equation}
where $\chi_1$ and $\chi_2$ are jointly Gaussian with $E(\chi
_1)^2=E(\chi_2)^2=\xi'(v)$ and $E\chi_1\chi_2= t\xi'(u)$
independent of $h$.
\end{proposition}

\begin{pf}
Using (\ref{ME0clem1eq3}), we have for $0<\delta<c$,
\begin{eqnarray*}
&&
\limsup_{n\rightarrow\infty}\alpha_n(0)
\\
&&\qquad\leq2\mathcal{P}(\xi,h)+\limsup_{n\rightarrow\infty}\mu_n \bigl( [0,c-
\delta] \bigr)\theta(1)+\bigl(\theta(v)-\theta(c-\delta)\bigr)_+
\\
&&\qquad=2\mathcal{P}(\xi,h)+\bigl(\theta(v)-\theta(c-\delta)\bigr)_+
\end{eqnarray*}
and this implies (\ref{ME0cprop3eq1}) by letting $\delta$ tend
to zero. For (\ref{ME0cprop3eq2}), we use (\ref{ME0clem1eq2}).
\end{pf}

Let us now turn to the study of some basic properties of $\varphi_c$.
Recall from (\ref{INTeq4}) and (\ref{ME0cprop3eq2}), for
fixed $0< v<1$, $\varphi_v$ is defined by
\[
\varphi_{v}(u,t)=\lim_{n\rightarrow\infty}\alpha_n'(0)
\]
for $0\leq u\leq v$ and $0\leq t\leq1$, where $\chi_1$ and $\chi_2$
are jointly Gaussian r.v.'s with
$E(\chi_1)^2=E(\chi_2)^2=\xi'(v)$ and $E\chi_1\chi_2=t\xi'(u)$
independent of $h$.
For given $k,\mathbf{m},\mathbf{q}$, let us recall the definition of
$(A_p)_{0\leq p\leq k+2}$ from (\ref{PReq1}). We also recall the
definitions of $(W_p)_{1\leq p\leq k+1}$ and $(\zeta_p)_{1\leq p\leq k+1}$
from Lemma~\ref{PRlem4}. Let us proceed with the following lemmas.

\begin{lemma}\label{ME0clem0}
Let $\varepsilon>0$ and $0<\delta<c$. Suppose that $l$ and $l'$ are
fixed integers with $1\leq l<l'\leq k+1$.
If $m_p\leq\varepsilon$ for every $1\leq p\leq l-1$, then
%
\begin{equation}
\label{ME0clem0eq1} E|W_1W_2\cdots
W_{l-1}-1|\leq M\varepsilon.
\end{equation}
If $c-\delta\leq q_p\leq q_{l'}$ for every $l\leq p\leq l'$, then
%
\begin{equation}
\label{ME0clem0eq2} EW_1W_2\cdots W_{l-1}|W_{l}W_{l+1}
\cdots W_{l'-1}-1|\leq M\sqrt{q_{l'}-c+\delta}.
\end{equation}
Here, $M$ depends only on $\xi$ and $h$.
\end{lemma}

\begin{pf}
Similar arguments as $(14.468)$ and $(14.469)$ in~\cite{Talag11} will
yield the announced results immediately.
\end{pf}

\begin{lemma}\label{ME0clem2} We have
%
\begin{eqnarray}
\label{ME0clem2eq1} E \biggl(\frac{\partial\Phi}{\partial
x}(h+\chi,c) \biggr)^2&=&c,
\\
\label{ME0clem2eq2} \xi''(c)E \biggl(
\frac{\partial^2 \Phi}{\partial x^2}(h+\chi,c) \biggr)^2&\leq&1,
\end{eqnarray}
where $\chi$ denotes a Gaussian r.v. with $E\chi^2=\xi'(c)$.
\end{lemma}

\begin{pf}
Recall that each $\mu_n$ corresponds to $k,\mathbf{m},\mathbf{q}$
and $\varepsilon$.
Since \mbox{$0<c<1$}, for each $n$ there exists some $0\leq s\leq k+1$ such that
$q_s\leq c<q_{s+1}$. Let us first claim that
%
\begin{equation}
\label{ME0clem2proofeq1} \lim_{n\rightarrow\infty}E|W_1\cdots
W_{s-1}-1|=0
\end{equation}
and if $\lim_{n\rightarrow\infty}q_{s+1}=c$, then we further have
%
\begin{equation}
\label{ME0clem2proofeq2} \lim_{n\rightarrow\infty}E|W_1\cdots
W_{s}-1|=0.
\end{equation}
Let $0<\delta<c$ be fixed.
Suppose that $1\leq l\leq s+1$ is the largest integer such that
$q_{l-1}\leq c-\delta$.
Since $\lim_{n\rightarrow\infty}\mu_n( [0,c-\delta
])=0$, we have that for large $n$,
$m_p\leq\varepsilon$ for every $0\leq p\leq l-1$. Using
(\ref{ME0clem0eq1}),
%
\begin{equation}
\label{ME0clem2proofeq3} E|W_1W_2\cdots
W_{l-1}-1|\leq M\varepsilon.
\end{equation}
On the other hand, since $c-\delta\leq q_p\leq c< q_{s+1}$ for $l\leq
p\leq s$, using (\ref{ME0clem0eq2}), we also get
%
\begin{equation}
\label{ME0clem2proofeq4}\quad EW_1W_2\cdots
W_{l-1}\llvert W_{l}W_{l+1}\cdots
W_{s-1}-1\rrvert\leq M\sqrt{q_s-c+\delta}\leq M\sqrt{
\delta}
\end{equation}
and
%
\begin{equation}
\label{ME0clem2proofeq5} EW_1W_2\cdots
W_{l-1}\llvert W_{l}W_{l+1}\cdots
W_{s}-1\rrvert\leq M\sqrt{q_{s+1}-c+\delta}.
\end{equation}
Using the triangle inequality, (\ref{ME0clem2proofeq3}) and
(\ref{ME0clem2proofeq4}), it follows that
\begin{eqnarray*}
&&
\limsup_{n\rightarrow\infty}E\llvert W_1W_2\cdots
W_{s-1}-1\rrvert
\\[-1pt]
&&\qquad\leq\limsup_{n\rightarrow\infty}EW_1W_2\cdots
W_{l-1} \llvert W_{l}W_{l+1}\cdots
W_{s-1}-1\rrvert
\\[-1pt]
&&\qquad\quad{} +\limsup_{n\rightarrow\infty}E|W_1W_2\cdots
W_{l-1}-1|
\\[-1pt]
&&\qquad\leq\lim_{n\rightarrow\infty}M\sqrt{\delta}+M\varepsilon
\\[-1pt]
&&\qquad=M\sqrt{\delta}.
\end{eqnarray*}
Similarly, if $\lim_{n\rightarrow\infty}q_{s+1}=c$, using the
triangle inequality, (\ref{ME0clem2proofeq3}) and
(\ref{ME0clem2proofeq5}),
we obtain
\begin{eqnarray*}
&&\limsup_{n\rightarrow\infty}E\llvert W_1W_2\cdots
W_{s}-1\rrvert
\\[-1pt]
&&\qquad\leq\limsup_{n\rightarrow\infty} EW_1W_2\cdots
W_{l-1}E_l\llvert W_{l}W_{l+1}\cdots
W_{s}-1\rrvert
\\[-1pt]
&&\qquad\quad{} +\limsup_{n\rightarrow\infty}E|W_1W_2\cdots
W_{l-1}-1|
\\[-1pt]
&&\qquad\leq\lim_{n\rightarrow\infty}M\sqrt{q_{s+1}-c+\delta}+M\varepsilon
\\[-1pt]
&&\qquad=\lim_{n\rightarrow\infty}M\sqrt{\delta}+M\varepsilon
\\[-1pt]
&&\qquad=M\sqrt{\delta}.
\end{eqnarray*}
Since $\delta>0$ is arbitrary, our claim follows.

Now, let us assume, without loss of generality, that the following
limits exist:
\[
\lim_{n\rightarrow\infty}q_s,\qquad \lim_{n\rightarrow\infty
}q_{s+1},\qquad
\lim_{n\rightarrow\infty}m_{s}
\]
and denote them by $c_-$, $c_+$ and $m_c$, respectively. If
$c_-<c<c_+$, then the first inequality implies $m_c=0$, which leads to
a contradiction since
the second inequality implies $m_c>0$. Thus, we may assume
%
\begin{equation}
\label{ME0clem10}
\mbox{either $c_-=c$}\quad \mbox{or}\quad \mbox{$c_+=c$.}
\end{equation}
Note that from the stationarity of $\mu_n$, $q_p=0$ if and only if
$p=0$, and also $q_{p}=1$ if and only if $p=k+2$.
If $q_s=0$ for all but finitely many $n$, then $s+1=1\leq k+1$ for
large $n$ and so $c_+=c$.
If $q_{s+1}=1$ for all but finitely many~$n$,\vadjust{\goodbreak} then $s=k+1$ for large
$n$ and so $c_-=c$. Finally, if $0<q_s$ and $q_{s+1}<1$ for infinitely
many $n$,
then these $s$ satisfy $1\leq s\leq k$ and (\ref{ME0clem10}).
Hence, in the following argument, we assume further that one of the
following cases holds:
\begin{longlist}[(iii)]
\item[(i)] $1\leq s\leq k+1$ for all $n$ and $c_-=c$.
\item[(ii)] $1\leq s+1\leq k+1$ for all $n$ and $c_+=c$.
\item[(iii)] $1\leq s\leq k$ for all $n$ and (\ref{ME0clem10}) holds.
\end{longlist}
If (i) holds, then from (\ref{PRlem4eq1}), (\ref{PRlem4eq2})
and (\ref{ME0clem2proofeq1}), we have
\begin{eqnarray*}
E \biggl(\frac{\partial\Phi}{\partial x}(h+\chi,c) \biggr)^2&=&
\lim_{n\rightarrow\infty}EA_s'(h+\chi_s)^2
\\
&=&\lim_{n\rightarrow\infty}EW_1\cdots W_{s-1}A_{s}'(h+
\chi_s)^2
\\
&=&\lim_{n\rightarrow\infty}q_s
\\
&=&c
\end{eqnarray*}
and
\begin{eqnarray*}
\xi''(c)E \biggl(\frac{\partial^2\Phi}{\partial x^2}(h+\chi,c)
\biggr)^2&=&\lim_{n\rightarrow\infty}\xi''(q_s)EA_s''(h+
\chi_s)^2
\\
&\leq&\limsup_{n\rightarrow\infty}\xi''(q_s)EW_1
\cdots W_{s-1}A_s''(h+
\chi_s)^2
\\
&\leq&1,
\end{eqnarray*}
where $\chi_s$ is Gaussian with $E(\chi_s)^2=\xi'(q_s)$.
If (ii) holds, again from (\ref{PRlem4eq1}) and
(\ref{PRlem4eq2}), we have
\begin{eqnarray*}
EW_1\cdots W_{s}A_{s+1}'(h+
\chi_{s+1})^2&=&q_{s+1},
\\
\xi''(q_{s+1})EW_1\cdots
W_{s}A_{s+1}''(h+
\chi_{s+1})^2&\leq&1+M\varepsilon^{1/6}.
\end{eqnarray*}
Using (\ref{ME0clem2proofeq2}) and proceeding as in (i), we
obtain the announced results, where $\chi_{s+1}$ is Gaussian with
$E(\chi_{s+1})^2=
\xi'(q_{s+1})$. Finally, for the case (iii), the same argument
completes our proof.
\end{pf}

\begin{proposition}\label{ME0clem3}
For each $0\leq t\leq1$, $\varphi_{v}(\cdot,t)$ is a convex function
on $ [0,v ]$.
For $0\leq u\leq c$ and $0\leq t\leq1$,
%
\begin{eqnarray}
\label{ME0clem3eq1} \frac{\partial\varphi_{c}}{\partial u}&\leq&0;
\\
\label{ME0clem3eq2} \frac{\partial\varphi_{c}}{\partial t}&\geq&\frac
{\xi'}{M},
\end{eqnarray}
where $M$ is a constant depending only on $\xi$ and $h$.
\end{proposition}

\begin{pf}
Define for each $n$,
\[
\varphi_{n,v}(u,t)=\alpha_n'(0)=E
\frac{\partial\Phi_n}{\partial
x}(h+\chi_1,v)\,\frac{\partial\Phi_n}{\partial x}(h+
\chi_2,v)-u\vadjust{\goodbreak}
\]
for $0\leq u\leq v$ and $0\leq t\leq1$, where $\chi_1$ and $\chi_2$
are jointly Gaussian with $E(\chi_1)^2=E(\chi_2)^2=\xi'(v)$
and $E\chi_1\chi_2=t\xi'(u)$ independent of $h$.
Again, without loss of
generality, we may assume that $v=q_a$ for some $1\leq a\leq k+1$.
Let $g,g_0^1,g_0^2,g_1^1,g_1^2$ be i.i.d. Gaussian r.v.'s with variance
$\xi'(q_a)$ such that for $i=1,2$,
\[
\chi_i= \bigl(g\sqrt{t}+g_0^i\sqrt{1-t}
\bigr)\sqrt{\frac{\xi'(u)}{\xi'(q_a)}}+g_1^i\sqrt{1-
\frac{\xi'(u)}{\xi'(q_a)}}.
\]
Then $\varphi_{n,v}(u,t)$ can be written as
\[
\varphi_{n,v}(u,t)=\phi_n \biggl(\frac{\xi'(u)}{\xi'(q_a)},t
\biggr)-u,
\]
where
\[
\phi_n(w,t)=EA_a' \bigl(V_1(w,t)
\bigr)A_a' \bigl(V_2(w,t) \bigr),
\]
and for $i=1,2$,
\[
V_i(w,t)=h+ \bigl(g\sqrt{t}+g_0^i\sqrt{1-t}
\bigr)\sqrt{w}+g_1^i\sqrt{1-w}.
\]
So for $0\leq u\leq v$ and $0\leq t\leq1$, by using Gaussian
integration by parts,
%
\begin{eqnarray}
\label{ME0cproofeq1} \frac{\partial\varphi_{n,v}}{\partial
u}(u,t)&=&t\xi''(u)
\Gamma_1(u,t)-1,
\\
\label{ME0cproofeq2} \frac{\partial^2\varphi_{n,v}}{\partial
u^2}(u,t)&=&t\xi^{(3)}(u)
\Gamma_1(u,t)+t^2\xi''(u)^2
\Gamma_2(u,t),
\\
\label{ME0cproofeq3} \frac{\partial\varphi_{n,v}}{\partial
t}(u,t)&=&\xi'(u)
\Gamma_1(u,t),
\\
\label{ME0cproofeq4} \frac{\partial^2\varphi_{n,v}}{\partial
u\,\partial
t}(u,t)&=&t\xi'(u)
\xi''(u)\Gamma_2(u,t),
\end{eqnarray}
where
\begin{eqnarray*}
\Gamma_1(u,t)&=&EA_a''(h+
\chi_1)A_a''(h+
\chi_2),
\\
\Gamma_2(u,t)&=&EA_a^{(3)}(h+
\chi_1)A_a^{(3)}(h+\chi_2).
\end{eqnarray*}
Since $A_a''>0$, we have $\Gamma_1>0$. Let us also observe that
\[
E \bigl(A_a^{(3)} \bigl(V_1(w,t) \bigr)|g,h
\bigr)=E \bigl(A_a^{(3)} \bigl(V_2(w,t)
\bigr)|g,h \bigr),
\]
which implies $\Gamma_2>0$. Thus, using these and from
(\ref{ME0cproofeq2}) and (\ref{ME0cproofeq4}), we obtain
%
\begin{eqnarray}
\label{ME0cproofeq5} \frac{\partial^2\varphi_{v}}{\partial u^2}&=&\lim
_{n\rightarrow
0}
\frac{\partial^2\varphi_{n,v}}{\partial u^2}\geq0,
\\
\label{ME0cproofeq6} \frac{\partial^2\varphi_{v}}{\partial u\,\partial
t}&=&\lim_{n\rightarrow\infty}
\frac{\partial^2\varphi_{n,v}}{\partial
u\,\partial t}\geq0.
\end{eqnarray}
Thus, the convexity of $\varphi_v(\cdot,t)$ follows from
(\ref{ME0cproofeq5}).
By (\ref{ME0clem2eq2}) and (\ref{ME0cproofeq1}),
we know $\frac{\partial\varphi_c}{\partial u}(c,1)\leq0$ and from
(\ref{ME0cproofeq6}), this implies $\frac{\partial\varphi
_c}{\partial u}(c,t)\leq0$.
So we obtain (\ref{ME0clem3eq1}) by using (\ref{ME0cproofeq5}).
Finally, (\ref{ME0clem3eq2}) can be easily obtained from
(\ref{PRlem3eq2}) and~(\ref{ME0cproofeq3}).
\end{pf}

\begin{pf*}{Proof of Proposition~\ref{ME0cprop1}}
Let $0\leq t<1$. Notice that if $u=0$, then $\chi_1$ and $\chi_2$ are
independent and from (\ref{PRlem3eq1}), it implies $\varphi_c(0,t)>0$.
Since $\varphi_c(c,1)=0$ by (\ref{ME0clem2eq1}) and $\frac
{\partial\varphi_c}{\partial t}(c,t)\geq\xi'(c)/M>0$ from
(\ref{ME0clem3eq2}),
we conclude that $\varphi_c(c,t)<0$ and so $\varphi_c(\cdot,t)$ has
a solution in
$ [0,c ]$. Suppose that $u_1,u_2$ with $0<u_1<u_2< c$ are
two solutions of $\varphi_c(\cdot,t)=0$ in $ [0,c ]$.
From Rolle's theorem, there exists some $u_3$ with $u_1<u_3<u_2$ such
that $\frac{\partial\varphi_c}{\partial u}(u_3,t)=0$.
Using the convexity of $\varphi_c(\cdot,t)$, it implies $\frac
{\partial\varphi_c}{\partial u}(u,t)\geq0$ for all $u_3\leq u\leq c$
and so $\varphi_c(c,t)\geq\varphi_c(u_2,t)=0$, which contradicts to
$\varphi_c(c,t)<0$.
\end{pf*}

\begin{pf*}{Proof of Proposition~\ref{ME0cprop2}}
Combining (\ref{INTeq4}), (\ref{PRcor1eq1}),
(\ref{ME0cprop3eq1}) and (\ref{ME0cprop3eq2}), we get
that for $u,v,t$ with $0\leq u\leq v<1$ and $0\leq t\leq1$,
%
\begin{equation}
\label{ME0cthm1proofeq1} p_{N,u}\leq2\mathcal{P}(\xi,h)-
\tfrac{1}{2}\varphi_v(u,t)^2+\bigl(\theta(v)-
\theta(c)\bigr)_+.
\end{equation}
Applying $v=c$ to this inequality, we obtain (\ref{ME0cthm1eq1}).
Suppose that $0\leq t<1$ is fixed.
It is easy to see that
$(u,v)\mapsto\varphi_v(u,t)$ is continuous on $0\leq u\leq v<1$.
Since $\varphi_c(c,t)<0$, there exists some $\gamma>0$
such that $\varphi_v(u,t)\leq\varphi_c(c,t)/2$ whenever $c\leq u\leq
v\leq c+\gamma$.
By the continuity of $\theta$, we may also let $\gamma$ be small
enough such that $\theta(v)-\theta(c)<\varphi_c(c,t)^2/16$
whenever $c\leq v\leq c+\gamma$. Therefore, we obtain
(\ref{ME0cthm1eq2}) from (\ref{ME0cthm1proofeq1}).
\end{pf*}
%

\section{\texorpdfstring{Proof of Proposition \protect\ref{MEposprop1}}{Proof of Proposition 5}}\label{sec6}

In this section our main goal is to establish an iterative inequality
that is used in the proofs of Propositions~\ref{MEposprop1} and
\ref{MEprop3}. Let us start by stating our main result as follows.
Suppose that $y_1$ and $y_2$ are jointly Gaussian r.v.'s
with $E(y_1)^2=E(y_2)^2=1$ and $Ey_1y_2=t\geq0$ independent of~$h$.
Define
\begin{eqnarray*}
F_{1}(x_1,x_2,w)&=&E \bigl(
\T(x_1+y_1\sqrt{w})-\T(x_2+y_2
\sqrt{w}) \bigr)^2,
\\
F_{-1}(x_1,x_2,w)&=&E \bigl(
\T(x_1+y_1\sqrt{w})+\T(x_2-y_2
\sqrt{w}) \bigr)^2
\end{eqnarray*}
for $x_1,x_2\in\Bbb{R}$ and $w\geq0$. For convenience, we sometimes
simply denote $F_1$ by~$F$.
Recall the constant $C$ stated in Lemma~\ref{PRlem1}. Set
$C_0=t(2(1+t)C^2)^{-1}$.
For $0<|u|\leq1$, let $\eta\in\{-1,+1 \}$ satisfy
$u=\eta|u|$. Then the following inequality holds.

\begin{proposition}\label{MEposprop2}
There exists a constant $K_1$ depending only on $C$ and $\xi$ such
that the following statement holds.
Suppose that $0<c_1<c_2<1$ and
%
\begin{equation}
\label{MEposprop2eq1} 0<\xi'(c_2)-
\xi'(c_1)<\min\biggl(\frac{1}{8},
\frac{1}{2(2C_0\xi'(1)+K_1)} \biggr);
\end{equation}
and $k,\mathbf{m},\mathbf{q}$ are such that for some $1\leq s\leq k+1$,
%
\begin{equation}
\label{MEposprop2eq2} q_s\leq c_1 \quad\mbox{and}\quad
m_s\geq\delta.
\end{equation}
Then we have
%
\begin{equation}
\label{MEposprop2eq3} p_{N,u}\leq2\mathcal{P}_k(
\mathbf{m},\mathbf{q})-C_0\delta K_2\int
_{c_1}^{c_2}EF_\eta\bigl(h,h,
\xi'(q)\bigr)\xi''(q)\,dq
\end{equation}
for every $u$ with $c_2\leq|u|\leq1$, where $K_2$ is a constant
depending only on $\xi$.
\end{proposition}

As consequences of Proposition~\ref{MEposprop2}, Propositions
\ref{MEposprop1} and~\ref{MEprop3} now follow.

\begin{pf*}{Proof of Proposition~\ref{MEposprop1}} Set $c_2=c'$. Let
us choose $c_1\in(c,c' )$ such that
(\ref{MEposprop2eq1}) holds and $\mu$ is continuous at $c_1$.
Since $c$ is the minimum of the support of $\mu$, $\mu(
[0,c_1 ])>0$.
From the definition of $\mu$, there exists a sequence of $\varepsilon
_n$-stationary measures $(\mu_n)$ such that $\mu_n\rightarrow\mu$
weakly and
$\mathcal{P}(\xi,h,\mu_n)\rightarrow\mathcal{P}(\xi,h)$. For each
$n$, $\mu_n$ corresponds to some $k,\mathbf{m},\mathbf{q}$.
We assume that $c_1$ is in the list of $\mathbf{q}$ and $c_1=q_s$ for
some $1\leq s\leq k+1$.
Then for large $n$,
\[
\mu_n\bigl( [0,q_{s} ]\bigr)=m_{s}\geq
\delta,
\]
where $\delta=\mu( [0,c_1 ])/2$.
We then apply Proposition~\ref{MEposprop2} to obtain for every
$c'\leq u\leq1$,
\[
p_{N,u}\leq2\mathcal{P}_k(\mathbf{m},\mathbf{q})-
\varepsilon^*,
\]
where $\varepsilon^*=C_0\delta K_2\int_{c_1}^{c_2}EF_1(h,h,\xi'(q))\xi''(q)\,dq$.
Since $0<t<1$, we have that $\varepsilon^*>0$. Letting $n$ tend to
infinity completes our proof.
\end{pf*}

\begin{pf*}{Proof of Proposition~\ref{MEprop3}}
Note that from the given condition, we have \mbox{$v\geq\delta$}. Let
$c_1=v-\delta$ and $c_2=v-\delta/2$.
Without loss of generality, we may assume that $\delta>0$ is small
enough such that
(\ref{MEposprop2eq1}) holds. Since (\ref{MEposprop2eq2}) is
satisfied and $|u|=v>c_2$, it follows that from (\ref{MEposprop2eq3}),
\[
p_{N,u}\leq2\mathcal{P}_k(\mathbf{m},\mathbf{q})-
\varepsilon^*(v)\leq2\mathcal{P}(\xi,h)-\bigl(\varepsilon^*(v)-2
\varepsilon_1\bigr)
\]
for $\varepsilon^*(v)=C_0\delta K_2\int_{v-\delta}^{v-\delta
/2}EF_{-1}(h,h,\xi'(q))\xi''(q)\,dq$. Clearly, $\varepsilon^*(\cdot)$ is
a continuous function on $ [\delta,1 ]$. Since $0<t\leq1$ and
$Eh^2\neq0$, $\varepsilon^*(v)>0$ for every $v\in[\delta,1 ]$. Thus,
$\min_{v\in[\delta,1 ]}\varepsilon^*(v)>0$ and the announced result
follows by letting $\varepsilon_1$ be sufficiently small.
\end{pf*}

At this moment, we explain the motivation of the proof of Proposition
\ref{MEposprop2}. Let us apply (\ref{GBasseq2}) to Theorem
\ref{GBthm2}
and recall the definitions of $(Y_p)_{0\leq p\leq k+2}$ and
$(y_p^1,y_p^2)_{0\leq p\leq k+1}$.
Using the independence of $y_p^1$ and $y_p^2$ for $\tau\leq p\leq k+1$
and decreasing
induction, one may clearly derive
\[
Y_\tau=A_\tau\biggl(h+\sum_{0\leq p<\tau}y_p^1
\biggr)+A_\tau\biggl(h+\sum_{0\leq p<\tau}y_p^2
\biggr).
\]
For $0\leq p< \tau$, from Lemma~\ref{GBlem1} and again using
decreasing induction, we also have
%
\begin{eqnarray}
\label{sec6eq1} Y_{p}&=&\frac{1+t}{m_p}\log E_p\exp
\frac{m_p}{1+t}Y_{p+1}
\nonumber
\\
&\leq&\frac{1+t}{m_p}\log E_p\exp\frac{m_p}{1+t}
\bigl(A_{p+1} \bigl(x_p^1+y_p^1
\bigr)+A_{p+1} \bigl(x_p^2+y_p^2
\bigr) \bigr)
\\
&\leq&\frac{1}{m_p}E_p\exp m_p A_{p+1}
\bigl(x_p^1+y_p^1 \bigr)+
\frac
{1}{m_p}E_p\exp m_p A_{p+1}
\bigl(x_p^2+y_p^2 \bigr),
\nonumber
\end{eqnarray}
where $x_p^j=h+\sum_{0\leq r<p-1}y_r^j$ for $j=1,2$. In particular, if
$p=0$, $Y_0\leq\break2EA_0(h)=2X_0$.
To prove (\ref{MEposprop2eq3}), we expect that when $0<t<1$,
equality will not hold in (\ref{sec6eq1}) and,
with the help of the condition (\ref{MEposprop2eq2}), the small
difference between the two sides will keep accumulating over $p$.
Let us emphasize that this should be true even in the absence of the
external field.
A~similar approach is also presented in Section 14.12 of Talagrand's
book~\cite{Talag11}, where he considered the case $t=1$ and
used the Cauchy--Schwarz inequality to quantify the difference. However,
in the case $0<t<1$, his argument no longer holds.
We then resort to another approach using the Gaussian interpolation technique.

Before we state our main estimate, for convenience, let us set up a definition.
Let $C_1>0$ be a constant and $y$ be a standard Gaussian r.v. Suppose
that $m$ and $\omega$ are two fixed
numbers with $0\leq m\leq1$ and $\omega\geq0$ and $A$ is a
real-valued function defined on $\Bbb{R}$ such that
\[
E\exp mA(x+y\sqrt{w}) \quad\mbox{and}\quad EA(x+y\sqrt{w})
\]
exist for $x\in\Bbb{R}$ and $0\leq w\leq\omega$.
We define
%
\begin{equation}
\label{sec6eq2} T(x,w)=\frac{1}{m}\log E\exp mA(x+y\sqrt{w}),
\end{equation}
where $y$ is standard Gaussian.
Here, if $m=0$, $T(x,w)$ is defined as $EA(x+y\sqrt{w})$.
Then we say that $A$ satisfies condition $\mathcal{A}(m,\w,C_1)$ if
%
\begin{eqnarray}
\label{sec6eq3} \biggl\llvert\frac{\partial T}{\partial x}\biggr
\rrvert&\leq&1,\qquad
\frac{1}{C_1\C^2x}\leq\frac{\partial^2 T}{\partial x^2}\leq\min\biggl
(1,\frac
{C_1}{\C^2 x}
\biggr),\nonumber\\[-8pt]\\[-8pt]
\biggl\llvert\frac{\partial^3 T}{\partial x^3}\biggr\rrvert
&\leq&4,\qquad \biggl\llvert
\frac{\partial^4T}{\partial x^4}\biggr\rrvert\leq8\nonumber
\end{eqnarray}
for all $x\in\Bbb{R}$ and $0\leq w\leq\w$.

\begin{proposition}\label{ImportantThm1} Suppose that $A$ satisfies
$\mathcal{A}(m,\w,C_1)$.
Let $y_1,y_2$ be jointly Gaussian r.v.'s with $Ey_1^2=Ey_2^2=1$ and
$Ey_1y_2=t\geq0$.
Let $K>0$ and $L\in\Bbb{N}$ be fixed constants. Suppose that $\alpha
_0,\alpha_1,\ldots,\alpha_{\ell}\geq0$.\vadjust{\goodbreak}
Then there exist constants $C_\ell^0,C_\ell^1,\ldots,C_\ell^{\ell
}$ satisfying
%
\begin{equation}
\label{ImportantThm1eq2} 0<C_\ell^0,C_\ell^1,\ldots,C_\ell^\ell\leq4\sum_{n=0}^\ell
\alpha_n+K_1
\end{equation}
for some constant $K_1$ depending only on $C_1$ and $L$
such that for any given numbers $x_1,x_2\in\Bbb{R}$, $0< m\leq1$,
$0\leq w\leq\min(1/8,\w,1/{2C_\ell^0} )$, $w_0=0$, and
$0\leq w_1,w_2,\ldots,w_\ell\leq L$,
the following inequality holds:
%
\begin{eqnarray}
\label{ImportantThm1eq1}
&&
\frac{1+t}{m}\log E\exp\frac{m}{1+t}
\Biggl(A(x_1+y_1\sqrt{w})+A(x_2+y_2
\sqrt{w})
\nonumber
\\
&&\hspace*{78pt}\qquad{} -\sum_{n=0}^\ell\alpha_n
F(x_1+y_1\sqrt{w},x_2+y_2
\sqrt{w},w_n)\Biggr)
\nonumber\\[-8pt]\\[-8pt]
&&\qquad\leq\sum_{j=1}^2\frac{1}{m}\log
E\exp mA(x_j+y_j\sqrt{w})
\nonumber
\\
&&\qquad\quad{} -\sum_{n=0}^\ell\biggl(
\alpha_n\bigl(1-wC_{\ell}^n\bigr)+
\frac
{C_0}{2}mw\delta_0(n) \biggr)F\bigl(x_1,x_2,
\bigl(1-\delta_0(n)\bigr)w+w_n\bigr),
\nonumber
\end{eqnarray}
where $C_0=t(2(1+t)C_1^2)^{-1}$ and we define $\delta_0(n)=1$ if $n=0$
and $0$ otherwise.
\end{proposition}

Let us explain how to use this inequality. Observe that the left-hand
side of (\ref{ImportantThm1eq1}) differs from
(\ref{GBlem1eq1}) by
the $\ell+1$ quantities
\[
\bigl(\alpha_nF(x_1,x_2,w_n)
\bigr)_{0\leq n\leq\ell}
\]
at the present stage.
Most of them will be preserved in the new stage by
\[
\bigl(\alpha_n\bigl(1-wC_\ell^n\bigr)F
\bigl(x_1,x_2,\bigl(1-\delta_0(n)
\bigr)w+w_n\bigr) \bigr)_{0\leq n\leq\ell}
\]
with the additional term
\[
\frac{C_0}{2}tmwF(x_1,x_2,0).
\]
So after one step, we obtain $(\ell+1)+1$ terms in the new stage.
Continued iterations of (\ref{ImportantThm1eq1}) lead to a sum of
these small quantities that
will converge to some positive number if $w$ is not too small at each iteration.
This is the main reason we need the growth control on $C_\ell^0,C_\ell
^1,\ldots,C_\ell^\ell$ through (\ref{ImportantThm1eq2}).

Now, we turn to the proof of Proposition~\ref{MEposprop2}. Let
$k,\mathbf{m},\mathbf{q}$ be a given triplet. Recall the definition
of $(A_p)_{0\leq p\leq k+2}$ from (\ref{PReq1}). We will need the
following lemma.

\begin{lemma}\label{sec6lem1}
For each $0\leq p\leq k+1$, $A_{p+1}$ satisfies $\mathcal{A}(m_p,\xi
'(q_{p+1}),C)$.
\end{lemma}

\begin{pf}
Let $0\leq p\leq k+1$ be fixed. Suppose that $0\leq w\leq\xi
'(q_{p+1})$. Note that $\xi'$ is strictly increasing on $
[0,\infty)$
and $\xi'(0)=0$. Let $q$ satisfy \mbox{$\xi'(q)=\xi'(q_{p+1})-w$}.
Set $k'=k+1-p$. Consider
\begin{eqnarray*}
&&\mathbf{m}'\mbox{: } m_0'=0,m_1'=m_p\mbox{,
and $m_n'=m_{n+p-1}$ for $2\leq n
\leq k'+1$},
\\
&&\mathbf{q}'\mbox{: } q_0'=0,q_1'=q\mbox{,
and $q_n'=q_{n+p-1}$ for $2\leq n
\leq k'+2$}.
\end{eqnarray*}
Let $(B_n)_{0\leq n\leq k'+2}$ be defined in the same way as
$(A_p)_{0\leq p\leq k+2}$ by using the triplet $k',\mathbf{m}',\mathbf{q}'$.
Then it should be clear that $B_2=A_p$ and so
\begin{eqnarray*}
B_1(x)&=&\frac{1}{m_1'}\log E\exp m_1'B_2
\bigl(x+y\sqrt{\xi'(q_{p+1})-\xi'(q)} \bigr)
\\
&=&\frac{1}{m_p}\log E\exp m_pA_p(x+y\sqrt{w}),
\end{eqnarray*}
where $y$ is a standard Gaussian r.v. Since $(B_n)_{0\leq n\leq k'+2}$
satisfies (\ref{PRlem1eq1}), this completes our proof.
\end{pf}

\begin{pf*}{Proof of Proposition~\ref{MEposprop2}}
Let $C$ be the constant in Lemma~\ref{PRlem1} and $L$ be the
smallest integer such that $L\geq\xi'(1)$.
Suppose that $K_1$ is obtained from Proposition~\ref{ImportantThm1}
by using $C_1=C$ and $L$.
Again, without loss of generality, we may assume that
$c_1=q_{s_1},c_2=q_{s_2}$, $u=q_a$ for $1\leq s_1<s_2\leq a\leq k+2$.
Moreover, for $s_1\leq p\leq s_2-1$,
%
\begin{equation}
\label{Ithm2eq5} 0<\xi'(q_{p+1})-\xi'(q_p)<
\frac{1}{2\gamma(s_2-s_1)}
\end{equation}
and for $0\leq p<s_1$,
%
\begin{equation}
\label{Ithm2eq11} 0<\xi'(q_{p+1})-\xi'(q_p)<
\frac{1}{2\gamma},
\end{equation}
where $\gamma:=\max(4,2C_0\xi'(1)+K_1)$.
Let us note that such $(q_p)_{0\leq p\leq k+2}$ exists by the
discussion right below Theorem~\ref{GBthm2}
and using the assumption (\ref{MEposprop2eq1}).
Let us consider the following numbers:
\begin{eqnarray*}
\lambda&=&0,
\\
\tau&=&a,
\\
\kappa&=&k+1,
\\
n_p&=&\frac{m_p}{t+1} \qquad\mbox{if $0\leq p<\tau$}\quad
\mbox{and}\quad
m_p\qquad\mbox{if $\tau\leq p\leq\kappa$},
\\
\rho_p&=&q_p \qquad\mbox{for $0\leq p\leq\kappa+1$}.
\end{eqnarray*}
From (\ref{GBasseq1}),
\[
p_{N,u}\leq2\log2+Y_0-\sum_{0\leq p\leq k+1}m_p
\bigl(\theta(q_{p+1})-\theta(q_p)\bigr).
\]
Recall the definition of $(y_p^1,y_p^2)_{0\leq p\leq\kappa}$ from
Theorem~\ref{GBthm2}.
We define $Y_a(x_1,x_2)=A_a(x_1)+A_a(x_2)$ and for $1\leq p<a$,
\[
Y_p(x_1,x_2)=\frac{1+t}{m_p}\log E\exp
\frac
{m_p}{1+t}Y_{p+1}\bigl(x_1+y_p^1,x_2+y_{p}^2
\bigr).
\]
Finally, set $Y_0(x_1,x_2)=EY_1(x_1+y_0^1,x_2+y_0^2)$. It is obvious
that from the definition $Y_0=EY_0(h,h)$.
From Proposition~\ref{ImportantThm1}, we know $Y_{s_2}(x_1,x_2)\leq
A_{s_2}(x_1)+A_{s_2}(x_2)$.
Set $\eta_p=\xi'(q_{p+1})-\xi'(q_p)$ for $0\leq p\leq k+1$.
We claim that for $s_1\leq p<s_2$,
%
\begin{equation}
\label{Ithm2eq1} Y_p(x_1,x_2)\leq
A_p(x_1)+A_p(x_2)-\sum
_{n=p}^{s_2-1}\beta_{n,p}F_\eta
\Biggl(x_1,x_2,\sum_{l=p}^{n-1}
\eta_l \Biggr),
\end{equation}
where
%
\begin{equation}
\label{addeq3} \beta_{n,p}=\frac{C_0}{2}m_n
\eta_n \biggl(1-\frac{1}{s_2-s_1} \biggr)^{n-p}.
\end{equation}
Here, we adapt the definition $\sum_{\ell=p}^{p'}u_\ell=0$ whenever
$p>p'$ that remains enforced thereafter. Let $s_1\leq p<s_2$ and
consider the following numbers:
%
\begin{eqnarray}
\label{Ithm2eq6} \ell&=&s_2-p-1,
\nonumber
\\
m&=&m_{p},
\nonumber\\[-8pt]\\[-8pt]
\alpha_0&=&0 \quad\mbox{and}\quad \alpha_n=\beta_{{n+p},p+1}
\qquad\mbox{for $1\leq n\leq\ell$},
\nonumber
\\
w_0&=&0 \quad\mbox{and}\quad w_n=\sum
_{l=p+1}^{n+p-1}\eta_l \qquad\mbox{for } 1\leq n\leq
\ell.
\nonumber
\end{eqnarray}
From the definition of $w_n$, we know that $0\leq w_n\leq\xi'(1)\leq
L$ for $0\leq n\leq\ell$.
Since $A_{p+1}$ satisfies $\mathcal{A}(m_p,\xi'(q_{p+1}),C_1)$,
applying (\ref{Ithm2eq6}) to Proposition~\ref{ImportantThm1},
we obtain $(C_\ell^n)_{0\leq n\leq\ell}$ that, from
(\ref{ImportantThm1eq2}), satisfies
%
\begin{equation}
\label{Ithm2eq10} C_\ell^0,C_\ell^1,\ldots,C_\ell^\ell\leq4\sum_{n=0}^\ell
\alpha_n+K_1\leq2C_0\sum
_{n=p+1}^{s_2-1}m_n\eta_n+K_1
\leq\gamma.
\end{equation}
Using (\ref{Ithm2eq5}) and (\ref{Ithm2eq10}), we know for
$0\leq n\leq\ell$,
%
\begin{equation}
\label{Ithm2eq15} C_\ell^n\eta_p\leq\gamma
\eta_p<\frac{1}{2(s_2-s_1)}<\frac{1}{s_2-s_1}.
\end{equation}
Take $w=\eta_p$. Notice that from (\ref{Ithm2eq5}) and
(\ref{Ithm2eq15}), $w\leq\min(1/8,\xi'(q_{p+1}),1/2C_\ell^0)$. If
$u>0$, then from (\ref{ImportantThm1eq1}),
(\ref{Ithm2eq6}) and (\ref{Ithm2eq15}), we obtain
(\ref{Ithm2eq1}) since
\begin{eqnarray*}
Y_p(x_1,x_2)&\leq& A_p(x_1)+A_p(x_2)
\\
&&{}-\frac{C_0}{2}mwF(x_1,x_2,0)-\sum
_{n=1}^\ell\alpha_n
\bigl(1-C_\ell^nw \bigr)F(x_1,x_2,w+w_n)
\\[-1pt]
&\leq& A_p(x_1)+A_p(x_2)-
\frac{C_0}{2}m_p\eta_pF(x_1,x_2,0)
\\[-1pt]
&&{} -\sum_{n=p+1}^{s_2-1}\beta_{n,p}F
\Biggl(x_1,x_2,\sum_{l=p}^{n-1}
\eta_l \Biggr).
\end{eqnarray*}
If $u<0$, then
%
\begin{eqnarray}
\label{Ithm2eq14} Ey_p^1\bigl(-y_p^2
\bigr)&=&t\bigl(\xi'(q_{p+1})-\xi'(q_p)
\bigr),
\nonumber
\\[-1pt]
A_{p+1}\bigl(x_2+y_p^2
\bigr)&=&A_{p+1}\bigl(-x_2-y_p^2
\bigr),
\\[-1pt]
F_{-1}\bigl(x_1+y_p^1,x_2+y_p^2,w_n
\bigr)&=&F \bigl(x_1+y_p^1,-x_2-y_p^2,w_n
\bigr)
\nonumber
\end{eqnarray}
and it follows by applying $(x_1,-x_2)$ instead of $(x_1,x_2)$ and
$(y_1,y_2)=(y_p^1,-y_p^2)$ to Proposition~\ref{ImportantThm1} that
\begin{eqnarray*}
Y_{p}(x_1,x_2)&\leq& A_p(x_1)+A_p(-x_2)
\\[-1pt]
&&{}-\frac{C_0}{2}mwF(x_1,-x_2,0)\\[-1pt]
&&{}-\sum
_{n=1}^\ell\alpha_n
\bigl(1-C_\ell^n w\bigr)F(x_1,-x_2,w+w_p)
\\[-1pt]
&\leq& A_p(x_1)+A_p(x_2)-
\frac{C_0}{2}m_p\eta_pF_{-1}(x_1,x_2,0)
\\[-1pt]
&&{} -\sum_{n=p+1}^{s_2-1}\beta_{n,p}F_{-1}
\Biggl(x_1,x_2,\sum_{l=p}^{n-1}
\eta_l \Biggr),
\end{eqnarray*}
where, again, we use (\ref{Ithm2eq15}) for the second inequality.
This completes the proof of our claim.

Next, we claim that for $0\leq p\leq s_1$,
%
\begin{eqnarray}
\label{Ithm2eq4} Y_p(x_1,x_2)&\leq&
A_p(x_1)+A_p(x_2)
\nonumber\\[-8pt]\\[-8pt]
&&{}-\exp\Biggl(-2\gamma\sum_{l=p}^{s_1-1}
\eta_l \Biggr)\sum_{n=s_1}^{s_2-1}
\beta_{n,s_1}F_\eta\Biggl(x_1,x_2,
\sum_{l=p}^{n-1}\eta_l \Biggr).
\nonumber
\end{eqnarray}
If $p=s_1$, then (\ref{Ithm2eq4}) holds by (\ref{Ithm2eq1}).
Suppose $0\leq p<s_1$.
Let us consider the following numbers:
%
\begin{eqnarray}
\label{Ithm2eq9} \ell&=&s_2-s_1,
\nonumber
\\
m&=&m_{p},
\nonumber\\[-8pt]\\[-8pt]
\alpha_0&=&0 \quad\mbox{and}\quad \alpha_n=\exp\Biggl(-2\gamma
\sum_{l=p+1}^{s_1-1}\eta_l \Biggr)
\beta_{{n+s_1-1},s_1}\qquad \mbox{for } 1\leq n\leq\ell,
\nonumber
\\
\quad\qquad w_0&=&0 \quad\mbox{and}\quad w_n=\sum
_{l=p+1}^{n+s_1-2}\eta_l\qquad \mbox{for }1\leq n\leq
\ell,
\nonumber
\end{eqnarray}
where $\beta_{n,p}$ is defined in (\ref{addeq3}). As in our first
claim, since $0\leq w_n\leq\xi'(1)\leq L$ for $0\leq n\leq k$ and
$A_{p+1}$ satisfies $\mathcal{A}(m_p,\xi'(q_{p+1}),C_1)$,
we can apply Proposition~\ref{ImportantThm1} using (\ref{Ithm2eq9})
to obtain $(C_\ell^n)_{n=0}^\ell$ that, from
(\ref{ImportantThm1eq2}), satisfies
%
\begin{equation}
\label{Ithm2eq12} C_\ell^0,C_\ell^1,\ldots,C_\ell^\ell\leq4\sum_{n=0}^\ell
\alpha_n+K_1\leq2C_0\sum
_{n=0}^{\ell}m_n\eta_n+K_1
\leq\gamma.
\end{equation}
We conclude from (\ref{Ithm2eq11}) and (\ref{Ithm2eq12}) that
%
\begin{equation}
\label{addeq6} C_\ell^n\eta_p\leq1/2
\end{equation}
for $0\leq n\leq\ell$.
Note that $1-x\geq\exp(-2x)$ if $x\leq1/2$. Using this,
(\ref{Ithm2eq12}), and (\ref{addeq6}) yield
%
\begin{equation}
\label{Ithm2eq13} 1-C_\ell^n\eta_p\geq\exp
\bigl(-2C_\ell^n\eta_p \bigr)\geq\exp(-2
\gamma\eta_p ).
\end{equation}
Set $w=\eta_p$. Notice that from (\ref{Ithm2eq11}) and
(\ref{addeq6}), $w\leq\min(1/8,\xi'(q_{p+1}),1/2C_\ell^0)$.
If $u>0$, using (\ref{ImportantThm1eq1}) and (\ref{Ithm2eq13}),
we obtain
\begin{eqnarray*}
Y_p(x_1,x_2)&\leq& A_p(x_1)+A_p(x_2)-
\frac{C_0}{2}m_p\eta_pF(x_1,x_2,0)
\\
&&{}-\exp\Biggl(-2\gamma\sum_{l=p+1}^{s_1-1}
\eta_l \Biggr)\sum_{n=s_1}^{s_2-1}
\beta_{n,s_1}\bigl(1-C_\ell^n\eta_p
\bigr)F \Biggl(x_1,x_2, \sum_{l=p}^{n-1}
\eta_l \Biggr)
\\
&\leq& A_p(x_1)+A_p(x_2)
\\
&&{}-\exp\Biggl(-2\gamma\sum_{l=p}^{s_1-1}
\eta_l \Biggr)\sum_{n=s_1}^{s_2-1}
\beta_{n,s_1}F \Biggl(x_1,x_2, \sum
_{l=p}^{n-1}\eta_l \Biggr).
\end{eqnarray*}
If $u<0$, we obtain (\ref{Ithm2eq4}) by using
(\ref{Ithm2eq14}), applying $(x_1,-x_2)$ instead of $(x_1,x_2)$ and
$(y_1,y_2)=(y_p^1,-y_p^2)$ to (\ref{Ithm2eq1}),
and a similar argument as in the case $u>0$. This completes the proof
of our second claim.

Now, let $p=0$ in (\ref{Ithm2eq4}) and note that $m_n\geq\delta
/2$ for $n\geq s_1$.
We then obtain
\begin{eqnarray*}
Y_0&=&EY_0(h,h)
\\
&\leq& 2EA_0(h)\\
&&{}-\frac{C_0\delta}{2}\exp\bigl(-2\gamma
\xi'(1) \bigr)
\\
&&\hspace*{11.5pt}{} \times\biggl(1-\frac{1}{s_2-s_1} \biggr)^{s_2-s_1}\sum
_{n=s_1}^{s_2-1}\bigl(\xi'(q_{n+1})-
\xi'(q_n)\bigr)EF_\eta\bigl(h,h,
\xi'(q_n) \bigr).
\end{eqnarray*}
Since we can partition $ [c_1,c_2 ]$ so that $\max_{s_1\leq
p\leq s_2-1}\eta_p$ is arbitrarily small,
by passing to the limit,
\begin{eqnarray*}
Y_0&\leq&2X_0-\frac{C_0\delta}{2}\exp\bigl(-2\gamma
\xi'(1)-1 \bigr)
\int_{c_1}^{c_2}EF_\eta\bigl(h,h,
\xi'(q)\bigr)\xi''(q)\,dq
\end{eqnarray*}
and we are done.
\end{pf*}

At the end of this section, we will prove Proposition
\ref{ImportantThm1} and we proceed by two lemmas.

\begin{lemma}\label{ImportantLem3}
For any $x_1,x_2\in\Bbb{R}$, $0\leq w\leq\frac{1}{8}$, and $w'\geq
0$, we have
%
\begin{equation}
\label{ImportantLem3eq1} F\bigl(x_1,x_2,w'+w
\bigr)\geq\tfrac{1}{2}F\bigl(x_1,x_2,w'
\bigr).
\end{equation}
\end{lemma}

\begin{pf}
First we prove that for $x_1,x_2\in\Bbb{R}$ and $0\leq w\leq1/4$,
%
\begin{equation}
\label{ImportantLem3proofeq1} F(x_1,x_2,w)
\geq(1-4w)F(x_1,x_2,0).
\end{equation}
If (\ref{ImportantLem3proofeq1}) holds, then
\[
F(x_1,x_2,w)\geq\tfrac{1}{2}F(x_1,x_2,0),
\]
whenever $x_1,x_2\in\Bbb{R}$ and $0\leq w\leq1/8$ and this implies
(\ref{ImportantLem3eq1}) since for $w'\geq0$,
\begin{eqnarray*}
F\bigl(x_1,x_2,w'+w\bigr)&=&EF
\bigl(x_1+y_1\sqrt{w'},x_2+y_2
\sqrt{w'},w\bigr)
\\
&\geq&\tfrac{1}{2}EF\bigl(x_1+y_1
\sqrt{w'},x_2+y_2\sqrt{w'},0
\bigr)
\\
&=&\tfrac{1}{2}F\bigl(x_1,x_2,w'
\bigr),
\end{eqnarray*}
where $y_1$ and $y_2$ are jointly Gaussian r.v.'s with
$E(y_1)^2=E(y_2)^2=1$ and $Ey_1y_2=t$.
To prove (\ref{ImportantLem3proofeq1}), for fixed $x_1,x_2$, let
us set $\varphi(w)=F(x_1,x_2,w)$.
Define $G(x,y)=(\T x-\T y)^2$. Using Gaussian integration by parts, we have
%
\begin{eqnarray}
\label{addeq10} \varphi'(0)&=&\tfrac{1}{2}
\bigl(G_{11}(x_1,x_2)+G_{22}(x_1,x_2)+2tG_{12}(x_1,x_2)
\bigr)
\nonumber
\\
&=&(\T x_1-\T x_2) \bigl(\T''
x_1-\T''x_2\bigr)+\bigl(
\T'x_1-\T'x_2
\bigr)^2
\nonumber\\[-8pt]\\[-8pt]
&&{} +2(1-t)\T'x_1\T'x_2
\nonumber
\\
&\geq& (\T x_1-\T x_2) \bigl(\T''
x_1-\T''x_2\bigr).
\nonumber
\end{eqnarray}
Since
\[
\T''x_1-\T''x_2=2(
\T x_1-\T x_2) \bigl((\T x_1+\T
x_2)^2-1-\T x_1\T x_2\bigr),
\]
using this equation together with (\ref{addeq10}) leads to
\begin{eqnarray*}
\varphi'(0)&\geq& 2(\T x_1-\T x_2)^2
\bigl((\T x_1+\T x_2)^2-1-\T x_1\T
x_2\bigr)
\\
&\geq& -4(\T x_1-\T x_2)^2.
\end{eqnarray*}
We may also compute the second derivative of $\varphi$ to see that
\[
\max_{0\leq w\leq1}\bigl|\varphi''(w)\bigr|/2\leq C,
\]
where $C$ is a constant independent of $t,w$$,x_1,x_2$.
So
\begin{eqnarray*}
F(x_1,x_2,w)&=&\varphi(w)
\\
&\geq&\varphi(0)+\varphi'(0)w-Cw^2
\\
&\geq& (1-4w)F(x_1,x_2,0)-Cw^2.
\end{eqnarray*}
Set $\delta_i=wi/N$. It is easy to see by induction
\[
F(x_1,x_2,\delta_i)\geq(1-4
\delta_1)^iF(x_1,x_2,0)-Ci
\delta_1^2
\]
for $1\leq i\leq N$. In particular, if we put $i=N$ and let $N$ tend to
infinity, we obtain
$F(x_1,x_2,w)\geq\exp(-4w)F(x_1,x_2,0)\geq(1-4w)F(x_1,x_2,0)$ and
this completes the proof.
\end{pf}

\begin{lemma}\label{ImportantLem2}
Suppose that $A$ is a function defined on $\Bbb{R}$ satisfying
\begin{eqnarray*}
\bigl|A'\bigr|&\leq&1,\qquad
\frac{1}{C_1\C x^2}\leq A''(x)
\leq\min\biggl(1,\frac{C_1}{\C^2 x} \biggr),\\
\bigl|A^{(3)}\bigr|&\leq&4,\qquad
\bigl|A^{(4)}\bigr|\leq8
\end{eqnarray*}
for some constant $C_1$. Let $y_1,y_2$ be jointly Gaussian r.v.'s with
$Ey_1^2=Ey_2^2=1$ and $Ey_1y_2=t\geq0$.
Let $K>0$ and $L\in\Bbb{N}$ be fixed constants. Suppose that $0\leq
\alpha_0,\alpha_1,\ldots,\alpha_\ell\leq K$.
Then there exist constants
\begin{eqnarray*}
&K_1 \mbox{ depending only on $C_1$ and $L$},&
\\
&\displaystyle  C_\ell^0,C_\ell^1,\ldots,C_\ell^\ell\leq\sum_{n=0}^\ell
\alpha_n+K_1&
\end{eqnarray*}
and
\[
C_\ell^{\ell+1} \mbox{ depending only on $\ell$ and $K$}
\]
such that for any given numbers $x_1,x_2\in\Bbb{R}$, $0<m\leq1$,
$0\leq w\leq1/8$, $w_0=0$, and $0\leq w_1,w_2,\ldots,w_\ell\leq L$,
the following inequality holds:
%
\begin{eqnarray}
\label{ImportantLem2eq1}
&&
\frac{1+t}{m}\log E\exp\frac{m}{1+t}
\Biggl(A(x_1+y_1\sqrt{w})+A(x_2+y_2
\sqrt{w})
\nonumber
\\
&&\hspace*{89pt}\quad{} -\sum_{n=0}^\ell\alpha_n
F(x_1+y_1\sqrt{w},x_2+y_2
\sqrt{w},w_n)\Biggr)
\nonumber\\
&&\qquad\leq\sum_{j=1}^2\frac{1}{m}\log
E\exp mA(x_j+y_j\sqrt{w})
\\
&&\qquad\quad{} -\sum_{n=0}^\ell\bigl(
\alpha_n \bigl(1-C_{\ell}^nw
\bigr)+C_0mw\delta_0(n) \bigr)F\bigl(x_1,x_2,
\bigl(1-\delta_0(n)\bigr)w+w_n\bigr)
\nonumber
\\
&&\qquad\quad{} +C_{\ell}^{\ell+1}w^2,
\nonumber
\end{eqnarray}
where $C_0=t(2(1+t)C_1^2)^{-1}$ and $\delta_0(n)=1$ if $n=0$ and $0$ otherwise.
\end{lemma}

\begin{pf}
The proof is based on the Gaussian interpolation technique. Suppose for
the moment that $(y_1,y_2)$ are jointly Gaussian with
$E(y_1)^2=E(y_2)^2\leq1/8$ and $Ey_1y_2=tE(y_1)^2$. Let $(z_1,z_2)$ be
an independent copy of $(y_1,y_2)$. Define $(z_1^0,z_2^0)=(0,0)$ and
for $1\leq n\leq\ell$, $(z_1^n,z_2^n)=(z_1,z_2)$. For convenience, we
set for $j=1,2$,
\begin{eqnarray*}
A_j(x)&=&A(x_j+x),
\\
\T_j(x)&=&\T(x_j+x),
\\
U_j(u)&=&y_j\sqrt{u}
\end{eqnarray*}
and for $j=1,2$ and $n=0,1,2,\ldots,\ell$,
\begin{eqnarray*}
V_{n,j}(u)&=&y_j\sqrt{u}+z_j^n
\sqrt{1-u},
\\
G_n(u)&=&F\bigl(x_1+V_{n,1}(u),x_2+V_{n,2}(u),w_n
\bigr),
\\
G_{n,j}(u)&=&F_j\bigl(x_1+V_{n,1}(u),x_2+V_{n,2}(u),w_n
\bigr),
\\
G_{n,ij}(u)&=&F_{ij}\bigl(x_1+V_{n,1}(u),x_2+V_{n,2}(u),w_n
\bigr),
\end{eqnarray*}
where $F_j$ is the partial derivative of $F$ with respect to the $j$th
variable and $F_{ij}$ means the second partial derivative of $F$
with respect to $i$th and then $j$th variables.
Define the interpolation functions
\begin{eqnarray*}
\varphi(u)&=&E_z\psi(u),
\\
\varphi_j(u)&=&\psi_j(u),\qquad j=1,2,
\end{eqnarray*}
where
\begin{eqnarray*}
\psi(u)&=&\frac{1+t}{m}\log E_y T(u),
\\
\psi_j(u)&=&\frac{1}{m}\log E_y
T_j(u)
\end{eqnarray*}
and
\begin{eqnarray*}
T(u)&=&\exp\frac{m}{1+t} \Biggl(A_1\bigl(U_1(u)
\bigr)+A_2\bigl(U_2(u)\bigr)-\sum
_{n=0}^\ell\alpha_n G_n(u)
\Biggr),
\\
T_j(u)&=&\exp m A_j\bigl(U_j(u)\bigr).
\end{eqnarray*}
Then
\begin{eqnarray*}
\varphi(1)&=&\frac{1+t}{m}\log E\exp\frac{m}{1+t} \Biggl(A_1(y_1)+A_2(y_2)
\\
&&\hspace*{99pt}{} -\sum_{n=0}^\ell\alpha_n
F(x_1+y_1,x_2+y_2,w_n)
\Biggr),
\\
\varphi(0)&=&\varphi_1(0)+\varphi_2(0)-\sum
_{n=0}^\ell\alpha_n E_zF
\bigl(x_1+z_1^n,x_2+z_2^n,w_n
\bigr).
\end{eqnarray*}
In the following, we will try to find an upper bound for $\varphi'(0)$.
Consider
\begin{eqnarray*}
\psi'(u)&=&\frac{1}{E_y T(u)}E_y \Biggl[\sum
_{j=1}^2 \Biggl(U_j'(u)A_j'
\bigl(U_j(u)\bigr)\\
&&\hspace*{77.5pt}{}-\sum_{n=0}^\ell
\alpha_nV_{n,j}'(u)G_{n,j}(u)
\Biggr)T(u) \Biggr]
\\
&=&\frac{1}{2}J_0(u)-\frac{1}{2}\sum
_{n=0}^\ell\alpha_n
\bigl(J_1^n(u)+J_2^n(u) \bigr),
\end{eqnarray*}
where
\[
J_0(u)=\frac{1}{E_yT(u)}E_y \biggl[\frac{1}{\sqrt{u}}
\bigl(y_1A_1'\bigl(U_1(u)
\bigr)+y_2A_2'\bigl(U_2(u)\bigr)
\bigr)T(u) \biggr]
\]
and for $j=1,2$ and $n=0,\ldots,\ell$,
\[
J_j^n(u)=\frac{1}{E_yT(u)}E_y \biggl[
\biggl(\frac{y_j}{\sqrt{u}}-\frac{z_j^n}{\sqrt{1-u}} \biggr
)G_{n,j}(u)T(u)
\biggr].
\]
Using Gaussian integration by parts on $y$, we have
\begin{eqnarray*}
J_0(0)
&=&E_y(y_1)^2\sum
_{j=1}^2 \Biggl(A_j''(0)+
\frac{m}{1+t}A_j'(0) \Biggl(A_j'(0)-
\sum_{n=0}^\ell\alpha_nG_{n,j}(0)
\Biggr) \Biggr)
\\
&&{} +\frac{tm}{1+t}E_y(y_1)^2
\Biggl(A_1'(0) \Biggl(A_2'(0)-
\sum_{n=0}^\ell\alpha_nG_{n,2}(0)
\Biggr)
\\
&&\hspace*{80.6pt}{} +A_2'(0) \Biggl(A_1'(0)-
\sum_{n=0}^\ell\alpha_nG_{n,1}(0)
\Biggr)\Biggr)
\\
&=&E_y(y_1)^2 \biggl(J_{0}^1(0)-
\frac{m}{1+t} \bigl(J_{0}^2(0)+J_{0}^3(0)
\bigr) \biggr),
\end{eqnarray*}
where
\begin{eqnarray*}
J_0^1(0)&=&\sum_{j=1}^2
\bigl(A_j''(0)+mA_j'(0)^2
\bigr)-\frac{mt}{1+t} \bigl(A_1'(0)-A_2'(0)
\bigr)^2,
\\
J_0^2(0)&=&\bigl(A_1'(0)-A_2'(0)
\bigr)\sum_{n=0}^\ell\alpha_n
\bigl(G_{n,1}(0)+tG_{n,2}(0)\bigr),
\\
J_0^3(0)&=&(1+t)A_2'(0)\sum
_{n=0}^\ell\alpha_n
\bigl(G_{n,1}(0)+G_{n,2}(0) \bigr).
\end{eqnarray*}
Let us try to find an upper bound for $J_0(0)$ first.
Since
\[
\frac{1}{C_1\C^2 x}\leq A''(x)\leq\frac{C_1}{\C^2 x},
\]
it is easy to see from (\ref{ImportantLem3eq1}) that
%
\begin{eqnarray}
\label{ImportantLem2proofeq15} \frac{1}{C_1^2}F(x_1,x_2,w_0)
&\leq& \bigl(A_1'(0)-A_2'(0)
\bigr)^2
\nonumber
\\
&\leq& C_1^2F(x_1,x_2,w_0)
\\
&\leq& 2^{8L+1}C_1^2E_zF
\bigl(x_1+z_1^n,x_2+z_2^n,w_n
\bigr).
\nonumber
\end{eqnarray}
Since
\begin{eqnarray*}
\frac{\partial}{\partial x}(\T_1x-\T_2y)^2&=&2\bigl(1-
\T_1^2x\bigr) (\T_1x-\T_2y),
\\
\frac{\partial}{\partial y}(\T_1x-\T_2y)^2&=&2\bigl(1-
\T_2^2y\bigr) (\T_2y-\T_1x)
\end{eqnarray*}
from the Cauchy--Schwarz inequality, we have
%
\begin{eqnarray}
\label{ImportantLem2proofeq4} E_zG_{n,j}(0)^2
&=&E_zF_j\bigl(x_1+z_1^n,x_2+z_2^n,w_n
\bigr)^2
\nonumber
\\
&=&E_z \bigl(2E_{y'}\bigl[ \bigl(1-\T_j^2
\bigl(z_j^n+y_j'\sqrt
{w_n}\bigr) \bigr)
\nonumber\\
&&\hspace*{39.1pt}{} \times\bigl(\T_1\bigl(z_1^n+y_1'
\sqrt{w_n}\bigr)-\T_2\bigl(z_2^n+y_2'
\sqrt{w_n}\bigr) \bigr)\bigr]\bigr)^2
\\
&\leq&4E_zE_{y'} \bigl(\T_1
\bigl(z_1^n+y_1'
\sqrt{w_n}\bigr)-\T_2\bigl(z_2^n+y_2'
\sqrt{w_n}\bigr) \bigr)^2
\nonumber
\\
&=&4E_zF\bigl(x_1+z_1^n,x_2+z_2^n,w_n
\bigr),
\nonumber
\end{eqnarray}
where $y_1',y_2'$ are jointly Gaussian r.v.'s with
$E(y_1')^2=E(y_2')^2=1$ and $Ey_1'y_2'=t$.
Straightforward computation yields
\[
\frac{\partial}{\partial x}(\T_1x-\T_2y)^2+
\frac{\partial
}{\partial y}(\T_1x-\T_2y)^2=-2(
\T_1 x-\T_2 y)^2(\T_1x+
\T_2 y)
\]
and this implies
%
\begin{equation}
\label{ImportantLem2proofeq17} E_z\bigl\llvert
G_{n,1}(0)+G_{n,2}(0)\bigr\rrvert\leq4E_zF
\bigl(x_1+z_1^n,x_2+z_2^n,w_n
\bigr).
\end{equation}
Now, combining (\ref{ImportantLem2proofeq15}),
(\ref{ImportantLem2proofeq4}), (\ref{ImportantLem2proofeq17}), and
using Jensen's inequality,
\begin{eqnarray*}
E_zJ_0^1(0)&\leq&\sum
_{n=1}^2\bigl(A_j''(0)+mA_j'(0)^2
\bigr)-\frac
{mt}{(1+t)C_1^2}F(x_1,x_2,w_0),
\\
E_z\bigl|J_0^2(0)\bigr|&\leq&\bigl\llvert
A_1'(0)-A_2'(0)\bigr\rrvert
\sum_{n=0}^\ell\alpha_n \bigl(
\bigl(E_zG_{n,1}(0)^2 \bigr)^{1/2}+t
\bigl(E_zG_{n,2}(0)^2 \bigr)^{1/2}
\bigr)
\\
&\leq&2(1+t)\bigl|A_1'(0)-A_2'(0)\bigr|
\sum_{n=0}^\ell\alpha_n
\bigl(E_zF\bigl(x_1+z_1^n,x_2+z_2^n,w_n
\bigr) \bigr)^{1/2}
\\
&\leq& 2^{4L+{3/2}}C_1(1+t)\sum_{n=0}^\ell
\alpha_nE_zF\bigl(x_1+z_1^n,x_2+z_2^n,w_n
\bigr)
\end{eqnarray*}
and
\[
E_z\bigl|J_0^3(0)\bigr|\leq4(1+t)\sum
_{n=0}^\ell\alpha_nE_zF
\bigl(x_1+z_1^n,x_2+z_2^n,w_n
\bigr).
\]
To sum up, we obtain that
%
\begin{eqnarray}
\label{ImportantLem2proofeq18}\quad E_zJ_0(0)
&\leq&
E_y(y_1)^2\sum_{j=1}^2
\bigl(A_j''(0)+mA_j'(0)^2
\bigr)
\nonumber
\\
&&{} +mE_y(y_1)^2\sum
_{n=0}^\ell\biggl(\alpha_n
\bigl(2^{4L+{3/2}}C_1+4\bigr)-\frac{t}{(1+t)C_1^2}
\delta_0(n) \biggr)
\\
&&\hspace*{10pt}{} \times E_zF\bigl(x_1+z_1^n,x_2+z_2^n,w_n
\bigr).
\nonumber
\end{eqnarray}

Next, let us turn to the computation of $J_1^n$. By using Gaussian
integration by parts on $y$, we obtain
\begin{eqnarray*}
&&
\frac{1}{\sqrt{u}}E_y \bigl(y_1G_{n,1}(u)T(u)
\bigr)
\\
&&\qquad=E_y(y_1)^2E_y
\bigl[G_{n,11}(u)T(u)\bigr]+E_y (y_1y_2
)E_y\bigl[G_{n,12}(u)T(u)\bigr]
\\
&&\qquad\quad{}+\frac{m}{1+t}E_y(y_1)^2E_y
\Biggl[G_{n,1}(u) \Biggl(A_1'
\bigl(U_1(u)\bigr)-\sum_{l=0}^\ell
\alpha_l G_{l,1}(u) \Biggr)T(u) \Biggr]
\\
&&\qquad\quad{}+\frac{m}{1+t}E_y (y_1y_2
)E_y \Biggl[G_{n,1}(u) \Biggl(A_2'
\bigl(U_2(u)\bigr)-\sum_{l=0}^\ell
\alpha_l G_{l,2}(u) \Biggr)T(u) \Biggr],
\end{eqnarray*}
and this implies that
%
\begin{eqnarray}
\label{ImportantLem2proofeq10}
&&\lim_{u\rightarrow0}E_z \biggl[
\frac{1}{\sqrt{u}E_yT(u)}E_y \bigl(y_1G_{n,1}(u)T(u)
\bigr) \biggr]
\nonumber\\[-8pt]\\[-8pt]
&&\qquad=E_y(y_1)^2 \biggl(E_{z}G_{n,11}(0)+tE_zG_{n,12}(0)+
\frac{m}{1+t} \bigl(I_1^n+tI_2^n
\bigr) \biggr),
\nonumber
\end{eqnarray}
where for $0\leq n\leq\ell$,
\begin{eqnarray*}
I_1^n&=&E_{z} \Biggl[G_{n,1}(0)
\Biggl(A_1'(0)-\sum_{l=1}^\ell
\alpha_l G_{l,1}(0) \Biggr) \Biggr],
\\
I_2^n&=&E_{z} \Biggl[G_{n,1}(0)
\Biggl(A_2'(0)-\sum_{l=1}^\ell
\alpha_l G_{l,2}(0) \Biggr) \Biggr].
\end{eqnarray*}
On the other hand, letting $u\rightarrow0$ and then using Gaussian
integration by parts on~$z$, we also have
%
\begin{eqnarray}
\label{ImportantLem2proofeq11}
&&
\lim_{u\rightarrow0}E_z \biggl[
\frac{1}{\sqrt{1-u}E_yT(u)}E_y \bigl(z_1^nG_{n,1}(u)T(u)
\bigr) \biggr]
\nonumber
\\
&&\qquad=E_z \biggl[\frac{1}{E_yT(0)}z_1^nG_{n,1}(0)E_yT(0)
\biggr]=E_z \bigl[z_1^nG_{n,1}(0)
\bigr]\\
&&\qquad=E_z\bigl(z_1^n\bigr)^2
\bigl(E_z G_{n,11}(0)+tE_zG_{n,12}(0)
\bigr).
\nonumber
\end{eqnarray}
So from (\ref{ImportantLem2proofeq10}) and
(\ref{ImportantLem2proofeq11}),
%
\begin{eqnarray}
\label{ImportantLem2proofeq12}
&&
\lim_{u\rightarrow0}E_zJ_1^n(u)
\nonumber
\\
&&\qquad=\lim_{u\rightarrow0}E_{y,z} \biggl[\frac{1}{E_yT(u)} \biggl(
\frac
{y_1}{\sqrt{u}}-\frac{z_1^n}{\sqrt{1-u}} \biggr)G_{n,1}(u)T(u) \biggr]
\\
&&\qquad=E_y(y_1)^2 \biggl(\delta_0(n)
\bigl(G_{0,11}(0)+tG_{0,12}(0) \bigr)+\frac{m}{1+t}
\bigl(I_1^n+tI_2^n \bigr)
\biggr).
\nonumber
\end{eqnarray}
We may also compute $\lim_{u\rightarrow0}E_zJ_2^n(u)$ and this yields
%
\begin{eqnarray}
\label{ImportantLem2proofeq2}
&&
\lim_{u\rightarrow0}E_zJ_1^n(u)+E_zJ_2^n(u)
\nonumber\\[-8pt]\\[-8pt]
&&\qquad=\delta_0(n)E_y(y_1)^2I_0+
\frac{mE_y(y_1)^2}{1+t} \bigl(\bigl(I_1^n+I_3^n
\bigr)+t\bigl(I_2^n+I_4^n\bigr)
\bigr),
\nonumber
\end{eqnarray}
where
$I_0=G_{0,11}(0)+G_{0,22}(0)+2tG_{0,12}(0)$ and
for $0\leq n\leq\ell$,
\begin{eqnarray*}
I_3^n&=&E_{z} \Biggl[G_{n,2}(0)
\Biggl(A_2'(0)-\sum_{l=1}^\ell
\alpha_l G_{l,2}(0) \Biggr) \Biggr],
\\
I_4^n&=&E_{z} \Biggl[G_{n,2}(0)
\Biggl(A_1'(0)-\sum_{l=1}^\ell
\alpha_l G_{l,1}(0) \Biggr) \Biggr].
\end{eqnarray*}
Let us now try to find a suitable lower bound for
(\ref{ImportantLem2proofeq2}).
Observe that
\begin{eqnarray*}
\T_1'(0),\T_2'(0)&\geq&0,
\\
\T_1'(0)-\T_2'(0)&=&-\bigl(
\T_1(0)-\T_2(0)\bigr) \bigl(\T_1(0)+
\T_2(0)\bigr),
\\
\T_1''(0)-\T_2''(0)&=&2
\bigl(\T_1(0)-\T_2(0)\bigr)
\\
&&{} \times\bigl(\bigl(\T_1(0)+\T_2(0)\bigr)^2-1-
\T_1(0)\T_2(0)\bigr).
\end{eqnarray*}
This implies
%
\begin{eqnarray}
\label{ImportantLem2proofeq7} I_0&=&2\bigl(\T_1(0)-
\T_2(0)\bigr) \bigl(\T_1''(0)-
\T_2''(0)\bigr)
\nonumber
\\
&&{} +2\bigl(\T_1'(0)-\T_2'(0)
\bigr)^2+4(1-t)\T_1'(0)\T_2'(0)
\nonumber
\\
&\geq&2\bigl(\T_1(0)-\T_2(0)\bigr) \bigl(
\T_1''(0)-\T_2''(0)
\bigr)
\\
&\geq&-4\bigl(\T_1(0)-
\T_2(0)\bigr)^2
\nonumber
\\
&=&-4F(x_1,x_2,w_0).
\nonumber
\end{eqnarray}
As for the upper bounds for $|I_1^n+I_3^n|$ and $|I_2^n+I_4^n|$, we write
\begin{eqnarray*}
I_1^n+I_3^n&=&I_{11}^n+I_{12}^n+I_{13}^n,
\\
I_2^n+I_4^n&=&I_{21}^n+I_{22}^n+I_{23}^n,
\end{eqnarray*}
where
\begin{eqnarray*}
I_{11}^n&=&E_{z} \bigl[G_{n,1}(0)
\bigl(A_1'(0)-A_2'(0) \bigr)
\bigr],
\\
I_{21}^n&=&E_{z} \bigl[G_{n,2}(0)
\bigl(A_1'(0)-A_2'(0) \bigr)
\bigr],
\\
I_{12}^n&=&E_{z} \bigl[A_2'(0)
\bigl(G_{n,1}(0)+G_{n,2}(0) \bigr) \bigr],
\\
I_{22}^n&=&E_{z} \bigl[A_2'(0)
\bigl(G_{n,1}(0)+G_{n,2}(0) \bigr) \bigr],
\\
I_{13}^n&=&-\sum_{l=0}^\ell
\alpha_l E_{z} \bigl[ \bigl(G_{n,1}(0)G_{l,1}(0)+G_{n,2}(0)G_{l,2}(0)
\bigr) \bigr],
\\
I_{23}^n&=&-\sum_{l=0}^\ell
\alpha_l E_{z} \bigl[ \bigl(G_{n,1}(0)G_{l,2}(0)+G_{n,2}(0)G_{l,1}(0)
\bigr) \bigr].
\end{eqnarray*}
Using (\ref{ImportantLem2proofeq15}),
(\ref{ImportantLem2proofeq4}) and Jensen's inequality, we have
%
\begin{eqnarray}
\label{ImportantLem2proofeq19}
\bigl|I_{j1}^n\bigr|&\leq&
\bigl(E_yG_{n,j}(0)^2 \bigr)^{1/2}\bigl|A_1'(0)-A_2'(0)\bigr|
\nonumber\\[-8pt]\\[-8pt]
&\leq&2^{4L+{3/2}}C_1E_zF\bigl(x_1+z_1^n,x_2+z_2^n,w_n
\bigr).
\nonumber
\end{eqnarray}
From (\ref{ImportantLem2proofeq17}), we also have
%
\begin{equation}
\label{ImportantLem2proofeq20} \bigl|I_{j2}^n\bigr|
\leq4E_zF\bigl(x_1+z_1^n,x_2+z_2^n,w_n
\bigr).
\end{equation}
To bound $|I_{j3}^n|$, we use $ab\leq(a^2+b^2)/2$ and
(\ref{ImportantLem2proofeq4}). Then this leads to
%
\begin{eqnarray}
\label{ImportantLem2proofeq21} \bigl|I_{j3}^n\bigr|&\leq&
\frac{1}{2}\sum_{l=0}^\ell
\alpha_l \bigl(G_{n,1}(0)^2+G_{n,2}(0)^2+G_{l,1}(0)^2+G_{l,2}(0)^2
\bigr)
\nonumber
\\
&\leq& 4 \Biggl(\sum_{l=0}^\ell
\alpha_l \Biggr)E_zF\bigl(x_1+z_1^n,x_2+z_2^n,w_n
\bigr)
\\
&&{} +4\sum_{l=0}^\ell\alpha_lE_zF
\bigl(x_1+z_1^l,x_2+z_2^l,w_l
\bigr).
\nonumber
\end{eqnarray}
Now, combining (\ref{ImportantLem2proofeq7}),
(\ref{ImportantLem2proofeq19}), (\ref{ImportantLem2proofeq20}) and
(\ref{ImportantLem2proofeq21})
together, we obtain from~(\ref{ImportantLem2proofeq2}),
%
\begin{eqnarray}
\label{ImportantLem2proofeq22}
&&
\sum_{n=0}^\ell
\alpha_n \bigl(E_zJ_1^n(0)+E_zJ_2^n(0)
\bigr)
\nonumber
\\
&&\qquad\geq-4\alpha_0\delta_0(n)E_y(y_1)^2F(x_1,x_2,w_0)
\nonumber\\[-8pt]\\[-8pt]
&&\qquad\quad{}-mE_y(y_1)^2\sum
_{n=0}^\ell\alpha_n \Biggl(8\sum
_{l=0}^\ell\alpha_l+2^{4L+{3/2}}C_1+4
\Biggr)\nonumber\\
&&\hspace*{73.5pt}\qquad\quad{}\times E_zF\bigl(x_1+z_1^n,x_2+z_2^n,w_n
\bigr).
\nonumber
\end{eqnarray}
From now on, we replace $(y_1,y_2)$ by $(y_1\sqrt{w},y_2\sqrt{w})$
with $E(y_1)^2=1$. Combining (\ref{ImportantLem2proofeq18})
and (\ref{ImportantLem2proofeq22}), we get
%
\begin{eqnarray}
\label{ImportantLem2proofeq23}\quad \varphi'(0) &\leq&
\frac{w}{2}\sum_{j=1}^2
\bigl(A_j''(0)+mA_j'(0)^2
\bigr)
\nonumber\\[-8pt]\\[-8pt]
&&{} +w\sum_{n=0}^\ell\bigl(
\alpha_nC_\ell^n-C_0m
\delta_0(n) \bigr)F\bigl(x_1,x_2,\bigl(1-
\delta_0(n)\bigr)w+w_n\bigr),
\nonumber
\end{eqnarray}
where $C_0=t(2(1+t)C_1^2)^{-1}$ and for $0\leq n\leq\ell$,
\[
C_\ell^n=4m\sum_{l=0}^\ell
\alpha_l+2^{4L+
{3/2}}mC_1+4m+2\delta_0(n).
\]
It is easy to compute that
\[
\varphi_j'(0)=\frac{w}{2}
\bigl(A_j''(0)+mA_j'(0)^2
\bigr).
\]
We may also use Gaussian integration by parts and the given conditions
on the first four derivatives
to compute the second derivatives of $\varphi_1$, $\varphi_2$ and
$\varphi$ and this yields
\[
\frac{1}{2}\max_{0\leq u\leq1} \bigl(\bigl|\varphi_1''(u)\bigr|+\bigl|
\varphi_2''(u)\bigr|+\bigl|\varphi''(u)\bigr|
\bigr)\leq C_\ell^{\ell+1}w^2,
\]
where $C_\ell^{\ell+1}$ depends only on $\ell$ and $K$.
Finally, we finish by using the mean value theorem and
(\ref{ImportantLem2proofeq23}),
\begin{eqnarray*}
\varphi(1)&\leq&\varphi(0)+\varphi'(0)+\frac{1}{2}
\max_{0\leq
u\leq1}\bigl|\varphi''(u)\bigr|
\\
&\leq&\varphi_1(1)+\varphi_2(1)
\\
&&{} -\sum_{n=0}^\ell\bigl(
\alpha_n\bigl(1-C_\ell^nw\bigr)+C_0mw
\delta_0(n) \bigr)F\bigl(x_1,x_2,\bigl(1-
\delta_0(n)\bigr)w+w_n\bigr)
\\
&&{} +C_\ell^{\ell+1}w^2.
\end{eqnarray*}
\upqed\end{pf}

\begin{pf*}{Proof of Proposition~\ref{ImportantThm1}}
Recall that Lemma~\ref{ImportantLem2} guarantees the existence of
constants $C_0,C_\ell^0,\ldots,C_\ell^\ell$, which satisfy
(\ref{ImportantThm1eq2}). From (\ref{ImportantLem2eq1}), we only need
to prove that $C_\ell^{\ell +1}$ can be eliminated. To do this, let
$\alpha_0,\ldots,\alpha_\ell\geq0$ and let $C_\ell^{\ell +1}$ be
obtained by using $K=C_0\omega+\max(\alpha_0,\ldots,\alpha_\ell)$ in
Lemma~\ref{ImportantLem2}. Let us keep $0<m\leq1$, $t\geq0$, $0\leq
w\leq\min(1/8,\omega,1/2C_\ell^0)$, $w_0=0$, and $0\leq
w_1,\ldots,w_\ell\leq L$ fixed. We use $\varphi(x_1,x_2,w)$ to denote
the left-hand side of (\ref{ImportantThm1eq1}). Recall the definition
of $T(x,w)$ from (\ref{sec6eq2}) using $A$ and $m$. Set
$\delta_i=wi/N$ for $1\leq i\leq N$. We claim that for large $N$, the
following inequality holds:
%
\begin{eqnarray}
\label{ImportantLem2proofeq9} \quad\varphi(x_1,x_2,
\delta_i)&\leq& T(x_1,\delta_i)+T(x_2,
\delta_i)
\nonumber\\[-8pt]\\[-8pt]
&&{} -\sum_{n=0}^\ell\beta_{n,i}F
\bigl(x_1,x_2,\bigl(1-\delta_0(n)\bigr)
\delta_i+w_n\bigr)+iC_\ell^{\ell+1}
\delta_1^2
\nonumber
\end{eqnarray}
for all $1\leq i\leq N$, where
\[
\beta_{n,i}=\delta_0(n)C_0m
\delta_1\sum_{j=1}^i
\bigl(1-C_\ell^0\delta_1 \bigr)^{j-1}+
\alpha_n\bigl(1-C_\ell^n\delta_1
\bigr)^i.
\]
If $i=1$, then (\ref{ImportantLem2eq1}) implies
(\ref{ImportantLem2proofeq9}).
Suppose that (\ref{ImportantLem2proofeq9}) holds for some $i$ with
$1\leq i<N$.
Then by using the induction hypothesis,\vadjust{\eject}\vspace*{-12pt}
%
\begin{eqnarray}
\label{addeq4}
&&\varphi(x_1,x_2,\delta_{i+1})
\nonumber
\\
&&\qquad=\frac{1+t}{m}\log E\exp\frac{m}{1+t} \bigl(\varphi(x_1+y_{1}
\sqrt{\delta_1},x_2+y_{2}\sqrt{
\delta_1},\delta_i) \bigr)
\nonumber
\\
&&\qquad\leq\frac{1+t}{m}\log E\exp\frac{m}{1+t} \nonumber\\[-8pt]\\[-8pt]
&&\qquad\quad\hspace*{0pt}{}\times\Biggl(T(x_1+y_{1}
\sqrt{\delta_1},\delta_i)+ T(x_2+y_{2}
\sqrt{\delta_1},\delta_i)
\nonumber\\
&&\hspace*{18.2pt}\qquad\quad{} -\sum_{n=0}^\ell\beta_{n,i}F
\bigl(x_1+y_{1}\sqrt{\delta_1},x_2+y_{2}
\sqrt{\delta_1},\bigl(1-\delta_0(n)\bigr)
\delta_{i}+w_n\bigr) \Biggr)\nonumber\\
&&\qquad\quad{}+iC_\ell^{\ell+1}
\delta_1^2.
\nonumber
\end{eqnarray}
Observe that from the definition $\beta_{n,i}\leq C_0wi/N+\alpha_n\leq
K$ for large $N$ and $T(\cdot,\delta_i)$ satisfies $\mathcal
{A}(m,\delta_{N-i},C_1)$ since $0\leq w\leq\omega$. Also, notice
\[
\delta_1=\frac{w}{N}\leq\frac{1}{N}\min\biggl(
\frac{1}{8},\omega,\frac{1}{2C_\ell^0} \biggr)\leq\min\biggl(
\frac{1}{8},\delta_{N-i},\frac{1}{2C_\ell^0} \biggr).
\]
Applying (\ref{ImportantLem2eq1}) to (\ref{addeq4}), we obtain
\begin{eqnarray*}
&&\varphi(x_1,x_2,\delta_{i+1})
\\
&&\qquad\leq T(x_1,\delta_{i+1})+T(x_2,
\delta_{i+1})
\\
&&\qquad\quad{} -\sum_{n=0}^\ell\bigl(
\delta_0(n)C_0m\delta_1+\beta_{n,i}
\bigl(1-C_\ell^n\delta_1\bigr) \bigr)F
\bigl(x_1,x_2,\bigl(1-\delta_0(n)\bigr)
\delta_{i+1}+w_n\bigr)
\\
&&\qquad\quad{} +(i+1)C_\ell^{\ell+1}\delta_1^2.
\end{eqnarray*}
Since
\begin{eqnarray*}
&&
\delta_0(n)C_0m\delta_1+
\beta_{n,i}\bigl(1-C_\ell^n\delta_1
\bigr)
\\
&&\qquad=\delta_0(n)C_0m\delta_1+
\delta_0(n)C_0m\delta_1\sum
_{j=1}^i \bigl(1-C_\ell^0
\delta_1 \bigr)^{j}+ \alpha_n
\bigl(1-C_\ell^n\delta_1\bigr)^{i+1}
\\
&&\qquad=\delta_0(n)C_0m\delta_1\sum
_{j=1}^{i+1} \bigl(1-C_\ell^0
\delta_1 \bigr)^{j-1}+ \alpha_n
\bigl(1-C_\ell^n\delta_1\bigr)^{i+1}
\\
&&\qquad=\beta_{n,i+1},
\end{eqnarray*}
this completes the proof of our claim. Letting $i=N$ in
(\ref{ImportantLem2proofeq9}) and then $N\rightarrow\infty$, we obtain
that
%
\begin{eqnarray}
\label{addeq5} \varphi(x_1,x_2,w) &\leq&
T(x_1,w)+T(x_2,w)
\nonumber
\\
&&{} -\sum_{n=0}^\ell\bigl(
\delta_0(n)C_0 mw\exp\bigl(-C_\ell^0w
\bigr)+\alpha_n\exp\bigl(-C_\ell^nw\bigr)
\bigr)
\\
&&\hspace*{26pt}{} \times F\bigl(x_1,x_2,\bigl(1-\delta_0(n)
\bigr)w+w_n\bigr).
\nonumber
\end{eqnarray}
Since $\exp(-C_\ell^0w)\geq1/\sqrt{e}\geq1/2$ for $0\leq w\leq
1/2C_\ell^0$ and also $\exp(-C_\ell^n w)\geq1-C_\ell^nw$ using
$\exp(-x)\geq1-x$ for $x\geq0$, plugging these results inside
(\ref{addeq5}), we are done.
\end{pf*}

\section*{Acknowledgements}

The author would like to thank Michel Talagrand for sharing a
preliminary version of his book~\cite{Talag11}, which motivated the
present paper. Also, he would like to thank an anonymous referee and
Associate Editor for giving several valuable comments regarding the
presentation of the paper.



\printaddresses

\end{document}